\documentclass[final]{siamltex}
\pdfoutput=1

\usepackage{amssymb}
\usepackage{amsmath}
\usepackage{graphicx}

\usepackage{subcaption}
\usepackage{textcomp}
\usepackage{algpseudocode}
\usepackage{algorithm}
\usepackage[labelformat=simple]{subcaption}
\usepackage{afterpage}

\algrenewcommand{\algorithmiccomment}[1]{ \{ #1 \} }

\newif\ifboldnumber

\algrenewcommand\alglinenumber[1]{%
\footnotesize\ifboldnumber\bfseries\fi\global\boldnumberfalse#1:}

\usepackage{amsmath}
\usepackage{amssymb}
\usepackage{amsfonts}
\usepackage{graphicx}
\usepackage{color}

\newtheorem{exercise}{Exercise}
\newtheorem{exampl}{Example}
\newtheorem{conjecture}[theorem]{Conjecture}
\newcommand{\comment}[1]{} % makes its argument disappear

\def\bq{\begin{quotation}}
\def\eq{\end{quotation}}

\definecolor{purple}{rgb}{0.5,0.0,0.5}
\definecolor{darkgreen}{rgb}{0.2,0.5,0.2}

\def\c{\gamma}

\def\D{\Delta}

\def\k{\kappa}

\def\O{\Omega}

\newcommand{\V}[1]{ \mathbf{#1} }    % vector

\newcommand{\Vo}{\V{0}}

\newcommand{\Va}{\V{a}}
\newcommand{\Vb}{\V{b}}
\newcommand{\Vc}{\V{c}}

\newcommand{\Ve}{\V{e}}
\newcommand{\Vf}{\V{f}}

\newcommand{\Vn}{\V{n}}
\newcommand{\Vp}{\V{p}}

\newcommand{\Vu}{\V{u}}

\newcommand{\Vx}{\V{x}}
\newcommand{\Vy}{\V{y}}
\newcommand{\Vz}{\V{z}}

\newcommand{\M}[1]{ \mathbf{#1} }  % matrix

\newcommand{\MA}{\M{A}}
\newcommand{\MB}{\M{B}}

\newcommand{\ME}{\M{E}}

\newcommand{\MI}{\M{I}}

\newcommand{\MK}{\M{K}}
\newcommand{\ML}{\M{L}}

\newcommand{\MN}{\M{N}}

\newcommand{\MP}{\M{P}}

\newcommand{\MRr}{\M{R}}

\newcommand{\MU}{\M{U}}
\newcommand{\MV}{\M{V}}
\newcommand{\MW}{\M{W}}
\newcommand{\MX}{\M{X}}

\newcommand{\cD}{{\mathcal D}}
\newcommand{\cS}{{\mathcal S}}

\newcommand{\cF}{{\mathcal F}} % AGM addition
\newcommand{\Rn}[1]{\mathbb{R}^#1}
\newcommand{\Rmn}[2]{\mathbb{R}^{#1 \times #2}}
\newcommand{\Cn}[1]{\mathbb{C}^#1}
\newcommand{\Cmn}[2]{\mathbb{C}^{#1 \times #2}}

\def\Rl{\mathbb{R}}

\def\nm#1{\|#1\|}
\def\2nm#1{\|#1\|_2}
\def\Fnm#1{\|#1\|_F}
\def\Sp#1{\mathrm{Span}(#1)}

\newcommand{\ars}[1]{\left[ \begin{array}{#1}}
\newcommand{\are}{\end{array} \right] }
\newcommand{\oars}[1]{\begin{array}{#1}}
\newcommand{\oare}{\end{array}}
\newcommand{\rars}[1]{\left( \begin{array}{#1}}
\newcommand{\rare}{\end{array} \right) }

\newcommand{\eqs}{\begin{eqnarray}}
\newcommand{\eqe}{\end{eqnarray}}
\newcommand{\eqsn}{\begin{eqnarray*}}
\newcommand{\eqen}{\end{eqnarray*}}

\def\defs{\begin{definition}}
\def\defe{\end{definition}}
\def\teos{\begin{theorem}}
\def\teoe{\end{theorem}}
\def\prfs{\begin{proof}}
\def\prfe{\end{proof}}
\def\exas{\begin{exampl}}
\def\exae{\end{exampl}}
\def\excs{\begin{exercise}}
\def\exce{\end{exercise}}
\def\cors{\begin{corollary}}
\def\core{\end{corollary}}
\def\cons{\begin{conjecture}}
\def\cone{\end{conjecture}}

\newcommand{\ens}{\begin{enumerate}}
\newcommand{\ene}{\end{enumerate}}

\newcommand{\its}{\begin{itemize}}
\newcommand{\ite}{\end{itemize}}

\newcommand{\des}{\begin{description}}
\newcommand{\dee}{\end{description}}

\def\wh{\widehat}

\def\wt{\widetilde}

\def\rtm{\textsuperscript{\textregistered}}

\title{Preconditioning Parametrized Linear Systems \thanks{This material is
based upon work supported
by the National Science Foundation under grant numbers
NSF-DMS 1025327 and NSF-DMS 1217156, and by the Air Force Office of Scientific Research under grant number AFOSR FA9550-12-1-0442}}

\author{Arielle Carr\thanks{
{Department of Mathematics, Virginia Tech, Blacksburg, VA.}
({\tt arielle5@vt.edu})}
\and
Eric de Sturler\thanks{Department of Mathematics, Virginia Tech, Blacksburg, VA.
({\tt sturler@vt.edu}).}
\and
Serkan Gugercin\thanks{Department of Mathematics  and Computational Modeling and Data Analytics Division, Academy of Integrated Science, Virginia Tech, Blacksburg, VA.
({\tt gugercin@vt.edu}).}}

\begin{document}

\maketitle

\begin{abstract}
Preconditioners are generally essential for fast convergence in the iterative solution
of linear systems of equations. However, the computation of a good preconditioner
can be expensive. So, while solving a sequence of many linear systems, it
is advantageous to recycle preconditioners, that is, update a previous
preconditioner and reuse the updated version. In this paper, we introduce a
simple and effective method for doing this.
We consider recycling preconditioners for both
the general case of sequences of linear systems
$\MA(\Vp_k) \Vx_k = \Vb_k$ as well as the
important special case
of the type $(s_k\ME + \MA)\Vx_k = \Vb_k$.
The right hand sides
may or may not change.

We update preconditioners by defining a map from a new matrix to a
previous matrix, for example, the first matrix in the sequence.  We then combine the preconditioner for this previous matrix with the map
to define the new preconditioner. This approach has several advantages.
{\em The update is independent from the original preconditioner, so
it can be applied to any preconditioner.}
The possibly high cost of an
initial preconditioner can be amortized over many linear solves.
The cost of updating the preconditioner is more or less constant and there is flexibility in balancing
the quality of the map with the computational cost.

In the numerical experiments section, we demonstrate good results for several
applications, in particular when using an algebraic multigrid preconditioner.
\end{abstract}



\pagestyle{myheadings} \thispagestyle{plain} \markboth{A. Carr, E. de Sturler,
and S. Gugercin}{Preconditioning Parametrized Linear Systems}

%{\color{blue} Text in blue will indicate changes made to the paper; highlighted for the reviewers' convenience.}

\section{Introduction}\label{sec:intro}
We discuss the efficient computation of preconditioners for sequences of
systems that change slowly. %\AGM{We don't reference the first equation again using the equation number - either remove number or reference it immediately below.}
We consider both the general case
\eqs\label{eq:SeqGenParamSys}
  \MA(\Vp_k)\Vx_k = \Vb_k ,
\eqe
as well as the important special case
\eqs\label{eq:SeqShiftSys}
  (s_k\ME + \MA)\Vx_k = \Vb_k ,
\eqe
where the right hand side(s) may or may not change.
The first class of matrices in (\ref{eq:SeqGenParamSys}) arises, for example, in topology optimization,
discussed later in this paper, where
the parameter vector $\Vp_k$ represents the changing densities in each element
(during the optimization). $\MA(\Vp_k)$ represents the finite element discretization of
a three-dimensional elasticity problem given a density distribution $\Vp_k$ \cite{bendsoe1994, BendSig-Bk_2003, rozvany1995,WangStu_2007}.   In addition,
we consider two sequences of linear systems
of the form (\ref{eq:SeqShiftSys}). One 
arises in model reduction, in particular, in the Iterative
Rational Krylov Algorithm (IRKA)  \cite{AnthBeat10,GugeAnth08,AntBG20} for finding the optimal shifts $s_k$; the other arises in a sequence of discretized 2D Helmholtz equations \cite{BorgGuge14,ErlaNabb08,ErlanVOos_04}.  Other applications where sequences of the form (\ref{eq:SeqShiftSys}) arise include
oscillatory and transient hydraulic tomography (OHT/THT)  \cite{CardBarr11}
and diffuse optical tomography (DOT) \cite{AghaKilm11,deStKilm11,KilmdeSt06,SaibBakh13}, though we do not consider these applications in the current paper.
For the second class of matrices, $s_k$ is a shift
(often related to a frequency), and the matrices $\MA$ and $\ME$ ($\ME \neq \MI$)
represent discretizations of partial differential equations, or more generally arise in the simulation of a dynamical system. In these applications, the matrices
$\MA$ and $\ME$ may also depend on a parameter vector $\Vp$,
but this is not considered here.  

For the special case of shifted systems where $\ME = \MI$, other approaches
for iterative solvers have been considered. Flexible preconditioning is used for problems of
this form in \cite{BaumVanG15,GuZhou07}.  In \cite{AhmaSzyl15} and \cite{KilmdeSt06}, the authors take
advantage of the shift invariance of Krylov subspaces.

Preconditioners are often essential for fast iterative solutions
of linear systems of equations, but the computation of a good preconditioner
can be expensive. Therefore, we consider {\em recycling preconditioners},
that is, updating a previous
preconditioner and reusing the updated version for solving a new linear
system.
For a sequence of linear systems, this may provide a substantial reduction
in cost compared with {\em recomputing} a new preconditioner for each system.  We can also periodically compute a new preconditioner from scratch, which includes the important case of solving all systems with a single preconditioner.  We refer to this as {\em reusing} the initial preconditioner.

The main idea for our approach comes from \cite{AhujClar11}.
Given a sequence of matrices, $\MA_k$, for $k = 0, 1, 2, \ldots$, and
a good preconditioner $\MP_0$ for $\MA_0$, such that
$\MA_0 \MP_0$ (or $\MP_0 \MA_0$) yields fast convergence, we could
compute for each system the ideal
map $\wh{\MN}_k$ such that
\eqs\label{eq:map}
  \MA_k \wh{\MN}_k = \MA_0
\eqe
and define the updated preconditioner as
\eqs\label{eq:updP0}
  \MP_k = \wh{\MN}_k\MP_0.
\eqe
Then, $\MA_0\MP_0 = \MA_1\MP_1 = \dots = \MA_k\MP_k$, and $\MA_k \wh{\MN}_k \MP_0 = \MA_0 \MP_0$ will yield
the same fast convergence for each $k$ as
the original preconditioned system.
In general, the matrix $\wh{\MN}_k \MP_0$ is never
computed; in an iterative method, we can multiply vectors
successively by these two matrices
(which does lead to some overhead).
If computing these maps
can be made cheap and the initial preconditioner is very good,
we obtain fast convergence for all systems at low cost. In some cases, as with the Flow matrices discussed in Section \ref{sec:ModRed}, the initial preconditioner may not result in fast convergence, for example, because the initial matrix $\MA_0$ is very ill-conditioned and/or far from 
diagonally dominant.  In such a case, the preconditioner for another, more appropriate, system matrix $\MA_j$, $j>0$, may be chosen, and we recycle $\MP_j$.

In this paper, we present a more general update scheme for recycling preconditioners
by mapping one matrix to another for which we have a good preconditioner.
This generalizes the approach in \cite{AhujClar11} (see next section) to any set of closely related
matrices. %\AGM{I include the next few sentences as a response to Review \# 1 who asked about the ``obvious choice" in computing and applying the map. This is my best interpretation of what they meant.}
%{\color{blue} Computing the ideal map as in (\ref{eq:map}), and subsequently applying it, would generally be too expensive.  In particular, 
%\eqs\label{eq:exactMap}
%\wh{\MN}_k = \MA_k^{-1}\MA_0
%\eqe
%involves the high cost of computing an exact inverse and (\ref{eq:exactMap}) will be dense. Then, the overhead associated with the matrix-vector products in an iterative method could be very high.}  
We do not seek an exact map, but rather compute an approximate map $\MN_k$ such that
\eqs\label{eq:approxMap}
  \MA_k \MN_k \approx \MA_0.
\eqe

%{\color{blue} For the special case of shifted systems where $\ME = \MI$, other approaches
%for iterative solvers have been considered. Flexible preconditioning is used for problems of
%this form in \cite{BaumVanG15,GuZhou07}.  In \cite{AhmaSzyl15} and \cite{KilmdeSt06}, the authors take
%advantage of the shift invariance of Krylov subspaces.}

In Section \ref{sec:SAMs}, we review previous work on updating preconditioners
and Sparse Approximate Inverses (SAI), the technique motivating our proposed update scheme.
We then introduce our update scheme, the Sparse Approximate Map (SAM) update.

In Section \ref{sec:theo}, we analyze sparsity patterns for SAMs.
Denser patterns can give more accurate maps, but they also increase the
cost to compute the map and to apply it (every iteration), which
needs to be compensated with a further reduction in iterations.
On the other hand, if effective maps can be found that are significantly sparser than the matrix,
recycling preconditioners will be highly favorable. We demonstrate this for 3D
elasticity problems arising in topology optimization.

In Section \ref{sec:impl}, we discuss efficient implementations of SAMs, as
well as an efficient MATLAB\rtm~ m-file implementation of the
ILUTP factorization \cite{Saad_ILUT94, Saad09}.  For our numerical experiments, we also use an algebraic multigrid (AMG) preconditioner implemented in MATLAB\rtm~\cite{HuLinZik2019}.\footnote{We thank Xiaozhe Hu for sharing his MATLAB\rtm~ m-file for the AMG preconditioner with us.}  %As we did not develop this code, a detailed discussion of its implementation is not included in this paper, but we provide our parameter choices in Section \ref{sec:topopt}.
%For the SAM updates, we provide an algorithm for preprocessing the sparsity
%pattern of the map and an algorithm for computing the map.
The computation of SAMs as well as multiplying by SAMs
is easily parallelized, though we do not address this in the current paper.

SAM updates are particularly effective for sequences of hard
problems where expensive preconditioners are needed for fast convergence and reusing a preconditioner is not effective.
This is the case for many KKT systems and problems
where an AMG preconditioner is needed
and the set-up phase is expensive.
Another example are matrix-free methods where, cost-wise, we can compute a matrix
only once to compute a preconditioner.
In this paper, we
demonstrate the effectiveness of SAMs for ILUTP-type
preconditioners \cite{Saad09,Saad03} 
and AMG preconditioners \cite{RugStu1985, RugStu1986, Stu1983, Stu2001}
\cite[Appendix A: An Introduction to Algebraic Multigrid, K. St\"{u}ben]{TroOosShu2001}. 
%StuAppend2000}
%\AGM{Trottenberg AND Stuben's appendix or just the latter?} 
ILUTP-type preconditioners are widely used,  and AMG preconditioners make a good case as they are very effective for hard problems but expensive to compute.
While we use ILUTP  and AMG preconditioners here to make the
case for recycling preconditioners, SAMs can be used with any
other preconditioner.

We note that
the purpose of this paper is not a time-wise comparison between
SAMs and any preconditioner; rather, they play complementary roles.  
Nevertheless, recycling a preconditioner does not make sense
if computing the SAM takes more time than computing a new preconditioner. So, time-wise comparisons between
computing the two are needed.
To make these comparisons
fair, we compare runtimes for (interpreted) m-files: our m-file for SAMs, our m-file implementation for ILUTP \cite{ilutp-m},  and  an m-file implementation for the  AMG preconditioner from \cite{HuLinZik2019}.
%\footnote{{\color{blue}We thank Xiaozhe Hu for sharing this MATLAB\rtm~ m-file for the AMG preconditioner with us.}}  
%% REMOVING DEFENSE OF ILUTP
Comparing with runtimes from
MATLAB\rtm's (compiled) {\tt ilu} (type `ilutp') has two important drawbacks.
First, compiled code runs much faster than interpreted code, which would
seriously skew the comparisons. Second, MATLAB\rtm's more recent
implementation of Saad's ILUTP determines
the amount of fill automatically, sometimes allowing large amounts of fill.  This makes
comparisons difficult and potentially makes
computing MATLAB\rtm's {\tt ilu} more expensive than necessary (see footnote 6 on page 15).% (\Cpageref{foot:ILU})). 

Recycling preconditioners by periodically updating an
initial or previous preconditioner with a SAM update can significantly
reduce total runtime compared with (1) computing a new preconditioner for
every system or periodically and (2) reusing a fixed preconditioner
for all systems.  If computing (and using) 
the SAM is cheaper than 
computing a new preconditioner and yields faster convergence,
recycling clearly wins in comparison (1).  This is not usually the case,
but it is possible; see Section \ref{sec:indef}, where
we demonstrate this for indefinite matrices arising from the Helmholtz
equation.  With respect to (1), generally the issue is whether the (typically) lower cost of computing the
SAMs outweighs the cost of additional iterations (due to a less
effective preconditioner) and the additional
matrix-vector product per iteration.
With respect to (2), the issue is whether the reduction in iterations
due to an improved preconditioner outweighs
the cost of computing the SAM update plus the extra cost per iteration.

In Section~\ref{sec:results}, we demonstrate the effectiveness of
SAM updates for applications from topology optimization and model reduction, as well as for indefinite matrices arising from Helmholtz equations, along with providing some
details of these applications. % Removed next section - seems out of place here. (Also, I think was in response to previous reviewer who had concerns about us comparing preconditioners and our update scheme
\comment{We show that recycling preconditioners by periodically updating an
initial or previous preconditioner with a SAM update can {\color{blue}significantly} 
reduce total runtime compared with (1) computing a new preconditioner for
every system or periodically and (2) reusing a fixed preconditioner
for all systems.  If updating the preconditioner is cheaper than
computing a new preconditioner and yields faster convergence,
recycling clearly wins in comparison (1).  This is not usually the case,
but it is possible; see Section \ref{sec:indef}, where
we briefly demonstrate this for indefinite matrices arising from Helmholtz
equations.  With respect to (1), the issue is whether the (generally) lower cost of computing the
SAMs outweighs the cost of additional iterations (due to a less
effective preconditioner) and the additional
matrix-vector product per iteration.
With respect to (2), the issue is whether the reduction in iterations
due to an improved preconditioner outweighs
the cost of computing the SAM update plus the extra cost per iteration.} 

% OUT OF PLACE HERE?  MOVED FURTHER UP IN THIS SECTION
% WHERE WE INTRODUCE THE FORM OF THE SYSTEMS
\comment{For the special case of shifted systems where $\ME = \MI$, other approaches
for iterative solvers have been considered. Flexible preconditioning is used for problems of
this form in \cite{BaumVanG15,GuZhou07}.  In \cite{AhmaSzyl15} and \cite{KilmdeSt06}, the authors take
advantage of the shift invariance of Krylov subspaces.}

Finally, in Section \ref{sec:concl}, we discuss conclusions and future work.

\section{Recycling Preconditioners}\label{sec:SAMs}
To avoid the potentially high cost of computing a
new preconditioner, we propose recycling an existing one using maps between matrices.
In \cite{AhujClar11}, this idea was exploited for a
Markov chain Monte Carlo (MCMC) process that resulted in a
long sequence of matrices changing by one row at a time.
So, $\MA_{k+1} = \MA_k + \Ve_{i_k}\Vu_k^T$,
where $i_k$ indicates which row changes.
%and $\Vu_k$ is the change in the row.
The ideal map for this case,
$\MI - (1+\Vu_k^T\MA_k^{-1}\Ve_{i_k})^{-1}\MA_k^{-1}\Ve_{i_k}\Vu_k^T$, comes for free,
as we already need to compute
$\Vu_k^T \MA^{-1}_k \Ve_{i_k}$ for the transition probability
in the MCMC process.  While this update is specific to the change in the matrix, the approach proposed in the present paper generalizes
the idea of recycling preconditioners to {\em any} set of closely related matrices.

Our preconditioner update is advantageous in several ways.
(1) To compute the map (ideal or approximate),
knowledge of the original preconditioner, $\MP_0$, is not required.
{\em Therefore, the map is independent of $\MP_0$ and can be applied
to any type of preconditioner.} (2) The cost of updating $\MP_0$ in this fashion
is more or less constant, and the potentially high cost of computing a good $\MP_0$ can be amortized over many linear solves.  (3) In practice, we do not
need the ideal map (\ref{eq:map}), but rather
an approximation, $\MN_k$, satisfying (\ref{eq:approxMap}).  We can balance the accuracy of $\MN_k$
with the cost of computing it.  \comment{Applying the map is cheaper than computing a new preconditioner from scratch and, in some cases, cheaper than reusing $\MP_0$.
We demonstrate this for several applications in Section~\ref{sec:results}.}

Our update scheme is motivated by the Sparse Approximate Inverse (SAI). So, we refer to it as a Sparse Approximate Map (SAM) update.  
The SAI is proposed in \cite{Bens73} and further developed in \cite{BensFred82, ChowSaad98,GrotHuck97, HollWath05, Huc1998} and references therein.  To define SAIs and SAMs we need the following definitions.
\defs \label{defn:pattern}
A sparsity pattern for $\Cmn nn$ is any subset of $\{1,2,\dots,n\}\times\{1,2,\dots,n\}$.
\defe

\defs \label{defn:subSpace}
Let $S$ be a sparsity pattern for $\Cmn nn$.  We define the subspace
$\cS \subseteq \Cmn nn$ as $\cS = \{\MX \in \Cmn nn$ $|$ $X_{ij} = 0$ if $(i,j) \not \in S\}$.
\defe

SAIs can be defined in several ways, but for the current discussion we use the following.
\defs
For $\MP, \MA \in \Cmn nn$, $\MI$ the identity matrix in $\Cmn nn$, and a given sparsity pattern $S$, the Sparse Approximate Inverse, $\MP$, for a matrix, $\MA$, is defined as the minimizer of
\eqs\label{eq:SAI}
\min_{\MP\in\cS}\Fnm {\MI-\M{AP}}.
\eqe
\defe

%\AGM{More SAI citations added here per Rev \#3's request. Blue text is in response to Rev \#2's request to put these numbers upfront.} 
The computation of a SAI (and variations, such as MSAI \cite{HuckKall_2007}) is easily parallelized as the computation of $n$ independent, and very small, least squares problems \cite{ AnzHucBraDon2018, GroHuc1995, GrotHuck97, HucKalRoySedWei2010, HuckKall_2007, Kall_2008}. 
For example, for the Flow application discussed in Section \ref{sec:ModRed}, the maximum size of the least squares problems to compute the SAM updates is $20 \times 7$, independent of the matrix dimension $n = 9~669$.
%as discussed in . 
SAM updates can analogously be computed in
parallel, which can be a substantial advantage on modern architectures.   However, we do not consider a parallel implementation of SAMs in this paper.  

Rather than considering the identity matrix in (\ref{eq:SAI}),
other work has focused on replacing it with another matrix, sometimes referred to as a target matrix \cite{HollWath05}.  The problem then becomes
\eqs\label{eq:targetSAI}
  \min_{\MP\in\cS}\Fnm {\MB-\M{AP}}.
\eqe
In \cite{ChowSaad98, HollWath05}, (\ref{eq:targetSAI}) is solved to improve a
preconditioner, $\MB^{-1}$, aiming to make $\M{APB}^{-1}$ closer to the identity matrix
than $\MA\MB^{-1}$.  As a preconditioner, $\MB^{-1}$ is assumed to be available, for example through an approximate factorization of $\MA$ (or $\MA^{-1}$).
However, the columns of $\MB$ must be computed in order to solve (\ref{eq:targetSAI}),
and the cost of constructing these columns can be relatively high \cite{ChowSaad98}.  Indeed, for a preconditioner like AMG, as the components (or factors) of an AMG preconditioner are not individually invertible, an {\em iterative solve for each column of $\MB$} is required, which is typically quite expensive (although one can use $\MA$ as a preconditioner). In order to reduce the cost of explicitly constructing $\MB$, iterative methods with numerical dropping are used to approximate the columns of $\MB$ in \cite{ChowSaad98}.
In special cases, the structure or type of the matrix can be exploited.
In \cite{HollWath05}, using the advection-diffusion equation and
targeting the Laplacian, Holland, et al. are able
to use a fast solver for the action of $\MB^{-1}$ with good results.
 
Our update scheme involves solving
\eqs\label{eq:introSAM}
  \MN_k & = & \arg \min_{\MN\in\cS}\Fnm {\MA_k\MN - \MA_0},
\eqe
and defining the updated preconditioner as
\eqs\label{eq:Pk}
\MP_k = \MN_k\MP_0.
\eqe
Here $\cS$ is the subspace defined by a chosen sparsity pattern $S$ as in Definition \ref{defn:subSpace}, and $\MA_0$ and $\MA_k$
are matrices from a given sequence and are relatively close to one another.  From (\ref{eq:introSAM}), we obtain the approximation (\ref{eq:approxMap}).  

%{\color{blue}Our technique does not determine the (approximate or exact) inverse of any matrix, but rather a map from one matrix to another.  This map can then be used to update {\it any} preconditioner type.} %Further, 
While the minimization in (\ref{eq:introSAM}) has a form similar to (\ref{eq:targetSAI}), there are fundamental differences in the approach to preconditioning. First, computing (\ref{eq:targetSAI}) involves improving an existing preconditioner, $\MB^{-1}$, for a fixed matrix, $\MA$, 
whereas our approach aims to compute a sequence of maps between nearby matrices to update a given preconditioner to an earlier matrix to compute  preconditioners for the new matrices; the maps depend only on the sequence of nearby matrices not on the existing preconditioner.
Second, for most preconditioners, $\Fnm{\MB-\MA}$ is quite large \cite{ChowSaad98}. So, typically
an accurate solution cannot be expected, unless the sparsity pattern of $\MP$ contains relatively many nonzeros. Of course, if an accurate
solution is obtained, the benefit is
faster convergence rather than maintaining the same convergence.  
On the other hand, our approach seeks to map one matrix to another
{\em closely related one} in a sequence of linear systems. So, $\Fnm{\MA_k-\MA_0}$ is likely to be relatively small, and therefore, typically, we expect a relatively accurate solution. Third, computing the columns of $\MB$ can be quite expensive.
For the preconditioners used in this paper, this is certainly the case.  When using an ILUTP preconditioner, computing $\MB$ as in (\ref{eq:targetSAI}) requires the relatively expensive product $\M{LU}$. Computing  the inverse of an AMG preconditioner requires an iterative solve per column.
%Even, computing the inverse of an AMG preconditioner is not just expensive, it may be (nearly) impossible to determine as the computations of the multigrid algorithm are generally not invertible (at a reasonable cost).
Our approach requires only the solution of the small least problems, as the columns of $\MA_0$, a previous matrix in the sequence of linear systems, are readily available.

%\AGM{The next paragraph, and derivations within it, are a response to Reviewer \# 2's suggestion to provide such a bound.} \EdS{We are responding to the reviewer, but it still needs to be a natural part of the paper. So, it needs a sentence, why we are discussion this.} 
Next, we consider some theoretical properties of the preconditioned systems using the map. These properties provide some insight in how well the approach may work as well as how to choose parameters. However, for practical use, the latter also requires an assessment of the cost of the map, and we leave this for future work.

First, we show that if $\MA_0 \MP_0$ is well-conditioned then a sufficiently good map, i.e., a sufficiently small residual $\MRr_k = \MA_k\MN_k-\MA_0$, implies that $\MA_k\MN_k \MP_0$ is well-conditioned as well. Note that ill-conditioning generally leads to poor convergence of iterative methods, and hence is important to avoid. For nonsymmetric systems the converse is unfortunately not true. For that reason, we also provide a discussion of the field of values for the preconditioned system using the map.  

Suppose we compute a map, $\MN_k$, such that
(for some submultiplicative norm $\|.\|$)
\eqs\label{eq:suffAssump}
\nm{\MA_k\MN_k\MP_0-\MA_0\MP_0} < \gamma\nm{(\MA_0\MP_0)^{-1}}^{-1},
\eqe
for a chosen $0 < \c < 1$.
Then 
\eqs \label{eq:boundnB-1}
  \|(\MA_k\MN_k\MP_0)^{-1}\| \leq \frac{\|\MI\|}{1-\c}\|(\MA_0\MP_0)^{-1}\| ,
\eqe
%and as $\|\MA_k\MN_k\MP_0\| \leq \|\MA_0\MP_0\| + \|\MA_k\MN_k\MP_0 - \MA_0\MP_0\| \leq
%\|\MA_0\MP_0\| + \c\|(\MA_0\MP_0)^{-1}\|^{-1}$,
and as $\|\MA_k\MN_k\MP_0\| = \|\MA_0\MP_0 + (\MA_k\MN_k\MP_0 - \MA_0\MP_0)\| \leq
\|\MA_0\MP_0\| + \c\|(\MA_0\MP_0)^{-1}\|^{-1}$,
\eqs \label{eq:boundKappa}
\k(\MA_k\MN_k\MP_0) & \leq & \frac{\|\MI\|}{1-\c} (\k(\MA_0\MP_0) + \c).
\eqe
The proof of (\ref{eq:boundnB-1}) is a minor variation of \cite[Corollary 7.19]{Dym13} or \cite[section 2.3.4]{golub1996matrix}.
The parameter $\gamma$ need not be small; for instance, taking $\gamma = \frac{1}{2}$ and using the matrix 2-norm, gives
\eqs\label{eq:boundKappa2}
\k_2(\MA_k\MN_k\MP_0) \leq (2\k_2(\MA_0\MP_0)+1).
\eqe
We do not control 
$\nm{\MA_k\MN_k\MP_0-\MA_0\MP_0} = \| \MRr_k \MP_0 \|$
directly, but could choose instead of (\ref{eq:suffAssump}) 
\eqs\label{eq:suffAssump2}
\| \MRr_k \| \leq \frac{\c}{\| \MP_0 \|} \|(\MA_0 \MP_0)^{-1} \|^{-1} .
\eqe

And so the bound in (\ref{eq:boundKappa2}) depends on controllable features: how well $\MP_0$ approximates $\MA_0^{-1}$, and the accuracy of $\MN_k$. An estimate of  $\|(\MA_0 \MP_0)^{-1} \|$ can be computed during the iterative solve with $\MA_0 \MP_0$. One obvious way to improve the update is to use a denser sparsity pattern when computing $\MN_k$.  However, the sparsity pattern of the SAM update must be chosen carefully in order to minimize runtime.  We examine choices in sparsity patterns in more detail in Section \ref{sec:theo}.

So, our approach can keep the condition number from increasing substantially. While this is important, it does not guarantee fast convergence for nonsymmetric systems. For this reason, we also consider a field of values related argument. If the field of values of a matrix is contained in a disk or an ellipse not containing the origin, convergence bounds for GMRES are available; see, e.g., \cite[pp. 56-57]{Gree97}.
If we take $\| \MRr_k \|_2 \leq \c / \| \MP_0 \|_2$ (note the use of the matrix 2-norm), we get for any unit vector $\Vy$,
\eqs \label{eq:bndFOV_01}
  \Vy^*(\MA_k\MN_k\MP_0)\Vy & = &
  \Vy^*(\MA_0\MP_0)\Vy + 
  \Vy^*(\MA_k\MN_k\MP_0 - \MA_0\MP_0)\Vy \quad \implies \\
\label{eq:inclFOV_01}  
  \Vy^*(\MA_k\MN_k\MP_0)\Vy & \in & 
      \cF(\MA_0\MP_0) + \cD(\Vo,\c) ,
\eqe
where $\cF$ denotes the field of values and
$\cD(\Vo,\c)$ is the closed disk centered at the origin with radius $\c$.
So, $\cF(\MA_k\MN_k\MP_0) \subseteq \cF(\MA_0\MP_0) + \cD(\Vo,\c)$. In this case $\c$ could be chosen based on an estimate of $\min |\cF(\MA_0\MP_0)|$ (obtained during the iterative solve with $\MA_0\MP_0$).

%While this discussion provides some insight into the presented results, practical use of these, typically pessimistic bounds, requires that we combine them with cost estimates. This will be future work.}

%how the spectral norm of the updated preconditioned system, $\MA_k\MN_k\MP_0$, compares with the that of the initial preconditioned system, $\MA_0\MP_0$.  In particular, letting $\MRr_k = \MA_k\MN_k-\MA_0$, if
%\eqs\label{eq:small}
%\nm{(\MA_0\MP_0)^{-1}\MRr_k\MP_0}\leq \gamma < 1,
%\eqe

This paper focuses on solving (\ref{eq:introSAM}) for each $k$
or selected $k$, but we can also
apply such a map incrementally. In that case, for the $k^{th}$ matrix we solve
\eqs \label{eq:incrementSAM2} % \label{eq:incrementSAM1}
  \MN_{k,j_m} & = & \arg \min_{\MN\in\cS} \|\MA_k\MN - \MA_{j_m}\|_F ,
\eqe
%and define
%\eqs
%  \MP_k & = & \MN_{k,j_m} \MP_{j_m},
%\eqe
and define  $\MP_k = \MN_{k,j_m} \MP_{j_m}$ with
$\MP_{j_m}  =  \MN_{j_m,j_l} \MP_{j_l},$ for $0 \leq j_l < j_m < k$ (and so on).
This includes the special case $j_m = k-1$, $j_l = j_m-1$ (and so on).

When preconditioning from the left,
we can take advantage of row-wise changes made to $\MA_k$, as is
the case with the QMC matrices described above.  We can define
\eqs \label{eq:LeftMap}
  \MN_k & = & \arg \min_{\MN\in\cS}\Fnm {\MN\MA_k - \MA_0},
\eqe
with $\MP_k = \MP_0\MN_k$.
In this case, the computation of the map can be made
significantly cheaper by considering only those rows of $\MA_k$ that
differ from $\MA_0$ when computing the least squares minimization.  The same applies when computing the map as in (\ref{eq:introSAM}) if only a few columns of the matrix change.  Using the maps as in (\ref{eq:incrementSAM2}) and (\ref{eq:LeftMap}) is future research.

%\AGM{Reviewer \#3 asked for more LU/IC references. I added a few more references on preconditioner updates for sequences of linear systems using the approach:  1) highlight a good update scheme for ILU (we already had this - and it was specifically noted by a reviewer from the first submission); 2) highlight one for AMG (this was the paper we discussed back in August); and 3) provide several references for update techniques for sequences of linear systems without as much discussion as (1) and (2).}

Other preconditioner update schemes for sequences of linear systems have been proposed. For these schemes, though, the initial preconditioner typically must be of a specific type.  A cheap update to the factorized approximate inverse (AINV) preconditioner is discussed in \cite{BenzBert03}.  Algorithms computing AINV preconditioners and AINV updates can be found in \cite{BellBert11,BenzBert03,BenzCull00,BenzHaws00,BenzKouh01,BenzTuma98,Rafi14}.  Several incremental, or iterative, update techniques to an ILU factorization
are described in \cite{CalgCheh10}. Two efficient update techniques for an ILUT decomposition in a matrix-free environment are described in \cite{TebTum2010}. The balanced incomplete factorization \cite{BruMarMasTum2008}, a modification of ILU, can be updated using the scheme proposed in \cite{CerMarMas2017}. For symmetric positive definite (SPD) matrices, we may consider an incomplete Choleksy (IC) factorization for our preconditioner.  An overview of efficient techniques for updating an IC preconditioner using low rank updates is provided in \cite{Berg2020}. For SPD matrices with a diagonal shift, we could use the update scheme to an $\M{LDL}^T$ factorization presented in \cite{BelSimSerBen2011}. For a sequence of complex symmetric systems defined by a diagonal shift, we can also consider the update to an $\M{LDL}^H$ factorization as presented in \cite{Ber2004}.

A cheap update for incomplete factorizations of sequences of linear systems is presented in \cite{AnztChow16}, which uses the iterative algorithm proposed in \cite{ChowPate}. In \cite{ChowPate}, the authors introduce a method for computing an ILU factorization in parallel and provide insight into the costs of computing an ILU on modern architectures.  The update scheme in \cite{AnztChow16} uses a factorization for a previous matrix in the sequence as an initial guess to the iterative algorithm from \cite{ChowPate}.  This results in an update to the previous factorization.  Good results are provided for a sequence of linear systems coming from model reduction \cite{reporedirect15}, however this update scheme requires that the initial preconditioner be an ILU factorization.  

%\AGM{New addition of the AMG preconditioner update scheme.} 
A scheme for recomputing an AMG preconditioner at a reduced cost for sequences of stochastic collocation systems is proposed in \cite{GorPow2012}.  Prior to its use as either a linear solver or a preconditioner, a (generally) expensive setup phase is necessary for AMG.  In \cite{GorPow2012}, the authors show that by reusing some of this setup information for a sequence of systems, the setup phase for the next system in the sequence can be performed faster than if doing so from scratch.  Again, this update scheme is specific to the AMG preconditioner.

%\AGM{Added the following to make clear that SAMs are not intended to compete with other, successful, update strategies for specific preconditioner types.} 
%{\color{blue} We stress that we do not aim to compete with any individual preconditioner-specific update scheme.  Rather, our objective is to propose a cheap update scheme that is agnostic to the preconditioner type, therefore providing flexibility in its application.}

%{\color{blue}A modified sparse approximate inverse (MSPAI) technique is proposed in \cite{HucKalRoySedWei2010} that takes advantage of the structure and sparsity patterns of matrices using probing techniques in order to improve a preconditioner, with good results shown for matrices arising in image deblurring (see also \cite{ChaMat1992} for an overview of probing techniques).}

\section{Experimental Analysis of Sparsity Patterns for SAMs}\label{sec:theo}
As the cost and effectiveness of a SAM depends on the sparsity pattern,
we analyze some choices here. In choosing sparsity patterns
for SAMs, we aim to balance the cost of computing
and applying the map with the number of GMRES iterations to reduce total runtime.
For SAIs, both adaptive and fixed sparsity patterns have been
considered. Computing the pattern and preconditioner adaptively tends to be expensive \cite{Benz02,BenzTuma99,Chow00,Chow01,Huck99}. Therefore,
we focus on fixed, a priori, sparsity patterns for SAMs, although
some previous work on SAIs has focused on making adaptive strategies
more efficient \cite{ChowSaad98,GrotHuck97,Huc2003}.
We examine choices in sparsity patterns using matrices from
two applications,  Flow and topology optimization (see Sections \ref{sec:ModRed} and \ref{sec:topopt}, respectively, for more detail on these  applications).

\begin{figure}[h!]
\captionsetup{font=scriptsize}
\centering
        \begin{subfigure}{.44\textwidth}
        \captionsetup{font=small}
        		\includegraphics[width=\linewidth]{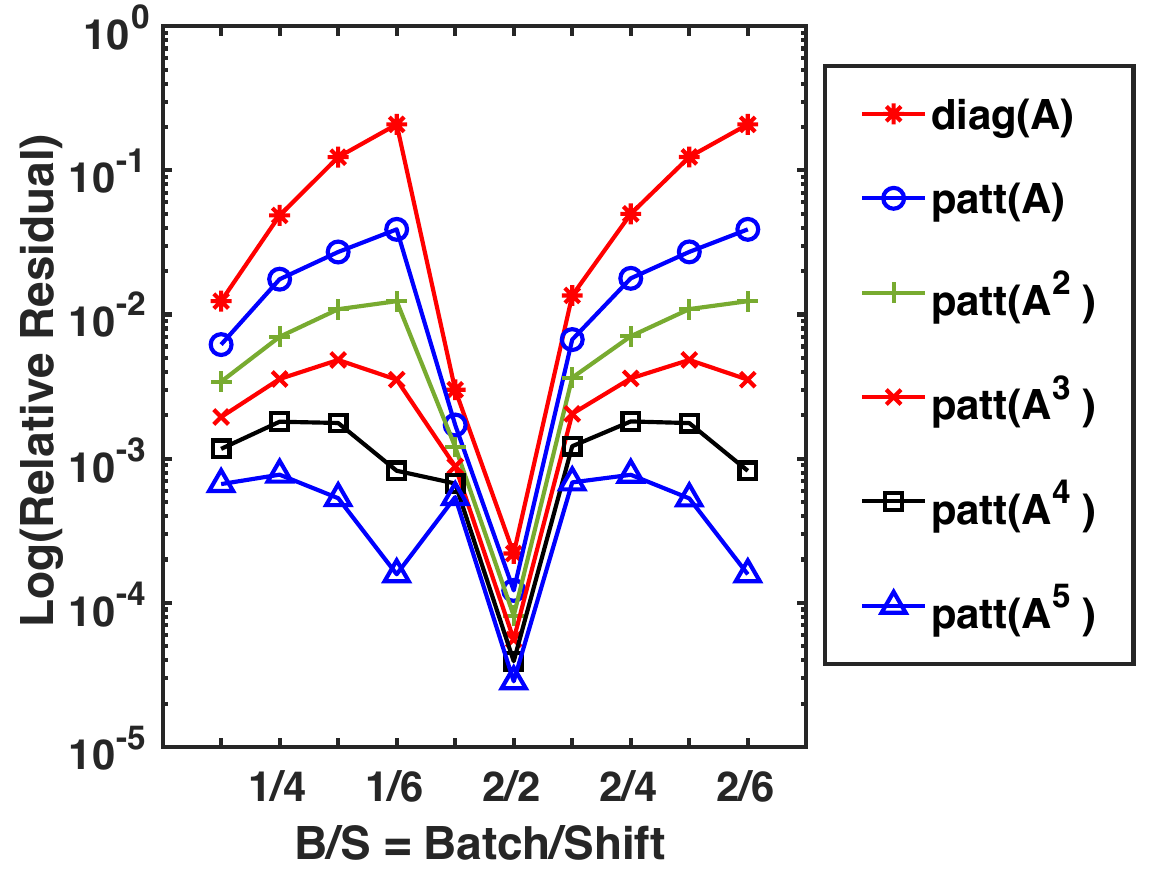}
		\subcaption{SAM relative residual norm vs shift for Flow matrices,
          for all patterns without sparsification.}
          \vspace{2mm}
		\label{fig:dense}
        \end{subfigure}
        \hspace{3mm} \begin{subfigure}{.42\textwidth}
        \captionsetup{font=small}
        		\includegraphics[width=\linewidth]{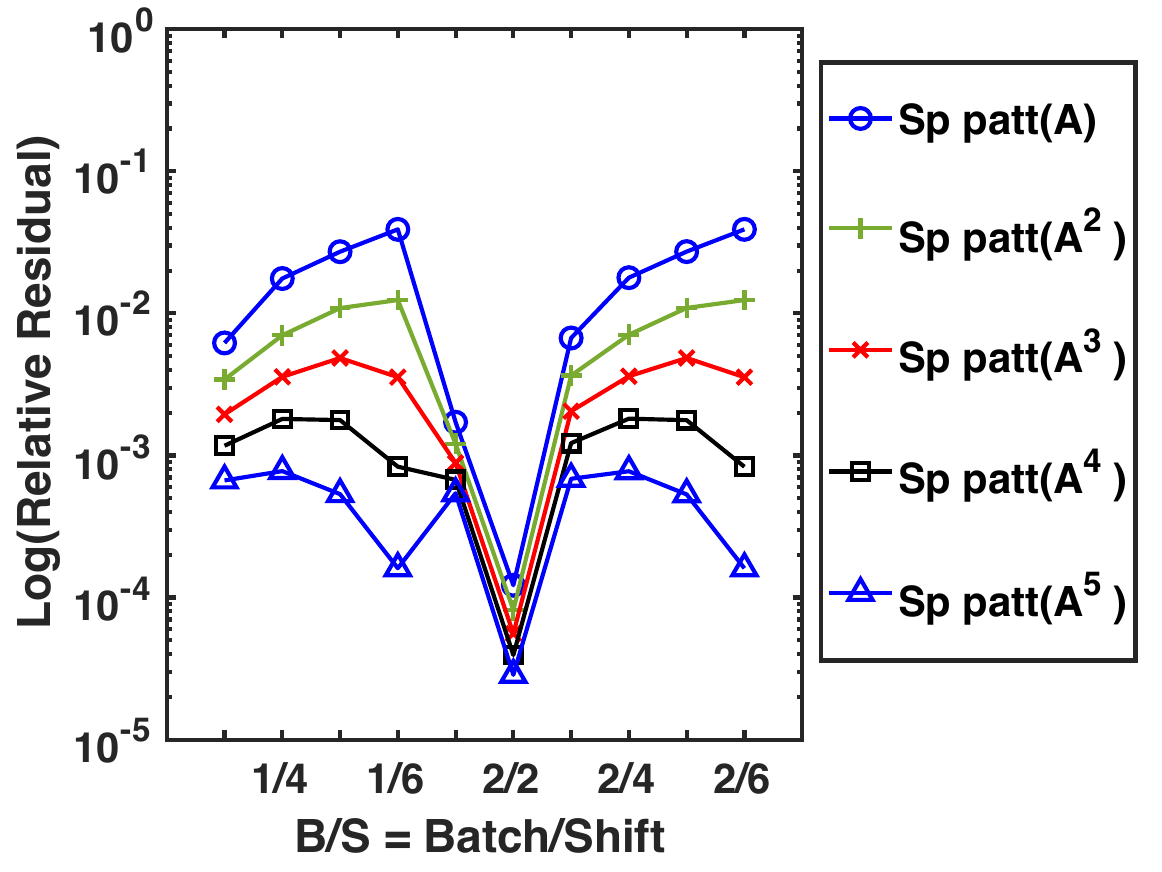}
		\subcaption{ SAM relative residual norm vs shift for Flow matrices,
          for all patterns with sparsification. ``Sp Patt" indicates ``Sparsified Pattern".}
		\label{fig:sparse}
        \end{subfigure}
        \vspace{-5mm}
\caption{Relative SAM residual, $\frac{\Fnm{\MA_k\MN_k-\MA_0}}{\Fnm{\MA_0}}$, of the
Flow matrices using a priori sparsity patterns for the SAM updates, $\MN_k$, for batch/shift 1/3 through 2/6. An ILUTP preconditioner is computed for batch/shift 1/2. Details of these linear systems and the notation are provided in Sections \ref{sec:theo} and \ref{sec:ModRed}.}
\vspace{-5mm}
\label{fig:Residuals}
\end{figure}

Sparsity patterns derived from powers of the matrix, $\MA$,
are a standard choice. This choice is based on the decay of the elements in the matrix
representing the discrete Green's function \cite{Chow00,Huck99,Tang99} associated
with the Laplace operator.
While the denser sparsity patterns of higher powers of $\MA$ may result in better maps,
solving (\ref{eq:introSAM}) becomes more costly.  We can alleviate this cost by
sparsifying the powers of $\MA$.  One possibility is to discard elements
of $\MA^k$ that are smaller in magnitude than a chosen threshold \cite{Chow00}. 
An alternative are sparsity patterns derived from the mesh on which the matrix
is based. %We follow this approach with topology optimization matrices based on a finite
%element discretization. 
Using the underlying mesh for the topology optimization application, we experiment with
sparsity patterns that are subsets of the sparsity pattern of the matrix, leading to maps that are much sparser than the system matrix.

We first analyze sparsity patterns defined by powers of the matrix, from $\MA$ to $\MA^5$,
sparsifications of those, and the diagonal pattern. The sparsifications are obtained using a global threshold
of $10^{-2}$.  To further minimize cost, we first sparsify $\MA$ using this global threshold, and then take powers of the logical matrix representing the nonzero pattern of this sparsified $\MA$.
We test these patterns
for the Flow matrices \cite{MoosRudn04,reporedirect15}. We study the relative residual norms of the resulting SAMs,
the time to compute the map, the number of GMRES iterations and GMRES runtime,
and the total solution time.  GMRES(200) is used for this application. As the sparsity pattern of the map becomes denser, we
expect the approximation in (\ref{eq:introSAM})
to become more accurate,
resulting in fewer GMRES iterations.

We consider eleven Flow matrices, $\MA_k = s_k\ME -\MA$, and corresponding
linear systems for a fixed right hand side $\Vb$ (see Section \ref{sec:ModRed}).  The Iterative Rational Krylov Algorithm (IRKA) \cite{GugeAnth08} computes optimal interpolation points, or shifts.  Each iteration of IRKA generates a batch of new shifts until IRKA converges.
%; so, many linear systems may need to be solved.  
We will refer to these matrices using the notation {\it batch/shift} (e.g., 1/2 is batch 1, shift 2).  We consider shifts 2 through 6 and 1 through 6 of batches 1 and 2, respectively (see Table \ref{table:FlowShifts} in Section \ref{sec:ModRed} for a list of all shifts). 

We compute an ILUTP preconditioner, $\MP_0$, for the first system in the
sequence, $\MA_0$, and we consider the powers of $\MA_0$ to derive
the denser patterns. For the Flow application, we compute and recycle the preconditioner for 1/2, not 1/1. For this first shift the matrix is quite ill-conditioned and not diagonally dominant, and this leads to a poor ILUTP preconditioner. So, we compute an ILUTP for this second matrix and recycle it for subsequent systems.
We define the residual of the map and its relative residual norm as
\eqsn
  \MRr_k = \MA_k\MN_k-\MA_0  \quad \mbox{and } \quad
  \frac{\Fnm{\MA_k\MN_k-\MA_0}}{\Fnm{\MA_0}} .
\eqen
Figure \ref{fig:Residuals} gives the relative residual norm for all shifts and all patterns,
and Tables \ref{table:Iters}  and \ref{table:ItersSparse} give,
for all patterns, GMRES iterations for selected shifts and (total) runtimes.

%\begin{center}
\begin{table}
%\parbox{0.99\linewidth}{
\captionsetup{font=scriptsize}
%\begin{center}
\centering
\scalebox{0.7}{
\begin{tabular}{|l|c|c|c|c|c|c|}
\hline
\textbf{Batch/Shift}&Diag&Patt $\MA$&Patt $\MA^2$&Patt $\MA^3$&Patt $\MA^4$&Patt $\MA^5$\\\hline
\textbf{1/2}&13&13&13&13&13&13\\\hline
\textbf{1/3}&20&19&18&17&16&15\\\hline
\textbf{1/4}&33&30&27&27&25&23\\\hline
\textbf{1/5}&108&53&49&44&38&31\\\hline
\textbf{1/6}&202&79&54&34&20&13\\\hline
\textbf{2/1}&496&488&1620&(5001)&(5001)&(5001)\\\hline
\textbf{2/2-2/6}&376&194&160&134&112&95\\\hline\hline
\multicolumn{7}{|c|}{{\bf ************ Totals ************}}\\\hline
%\textbf{2/2}&12&12&12&12&12&12\\\hline
%\textbf{2/3}&21&20&18&17&17&16\\\hline
%\textbf{2/4}&33&30&27&27&25&23\\\hline
%\textbf{2/5}&108&53&49&44&38&31\\\hline
%\textbf{2/6}&202&79&54&34&20&13\\\hline
\textbf{Total Iter}&1235&863&1928&5257&5212&5178\\\hline
\textbf{SAM Time (s)}&1.04&2.04&3.57&9.43&19.33&41.42\\\hline
\textbf{GMRES Time (s)} &3.16&1.47&3.49&11.70&14.41&16.36\\\hline
\textbf{Total Time (s) }&4.21&3.52&7.07&21.13&33.73&57.78\\\hline
\textbf{nnz($\MN_k$)/n}&1&6.97&18.94&37.02&61.39&92.36\\\hline\hline
\multicolumn{7}{|c|}{{\bf ********** Totals without 2/1 **********}}\\\hline
\textbf{Total Iter}&739&375&308&256&211&177\\\hline
\textbf{SAM Time (s)}&0.95&1.85&3.22&8.51&17.41&37.30\\\hline
\textbf{GMRES Time (s)} &2.44&0.75&0.60&0.56&0.54&0.53\\\hline
\textbf{Total Time (s) }&3.39&2.59&3.83&9.07&17.95&37.82\\\hline
\end{tabular}
}
\caption{GMRES iterations and total runtimes for selected shifts of the Flow matrices using a priori sparsity patterns for the SAM update.  An ILUTP preconditioner is computed for the system 1/2, with SAM updates computed for all remaining shifts 1/3-2/6. ``Patt" indicates ``"Pattern of".  ``SAM time" is the total amount of time spent computing all ten maps (but not including the time to compute the initial ILUTP). The initial ILUTP takes 0.8 s to compute.  Totals are given for the entire sequence, as well as for the sequence omitting system 2/1. (5001) indicates GMRES did not converge in the maximum allowed iterations.}
\label{table:Iters}
\vspace{-7mm}
\end{table}
%\end{center}
%
%\begin{center}
\begin{table}[!ht]
%\parbox{0.99\linewidth}{
\captionsetup{font=scriptsize}
%\begin{center}
\centering
\scalebox{0.7}{
\begin{tabular}{|l|c|c|c|c|c|}
\hline
\textbf{Batch/Shift}&Sp Patt $\MA$&Sp Patt $\MA^2$&Sp Patt $\MA^3$&Sp Patt $\MA^4$&Sp Patt $\MA^5$\\\hline
\textbf{1/2}&13&13&13&13&13\\\hline
\textbf{1/3}&19&18&17&16&15\\\hline
\textbf{1/4}&30&27&27&25&23\\\hline
\textbf{1/5}&53&50&45&40&34\\\hline
\textbf{1/6}&81&62&48&39&35\\\hline
\textbf{2/1}&(5001)&1437&(5001)&(5001)&(5001)\\\hline
\textbf{2/2-2/6}&196&170&149&133&120\\\hline\hline
%\textbf{2/2}&12&12&12&12&12\\\hline
%\textbf{2/3}&20&18&17&17&16\\\hline
%\textbf{2/4}&30&28&27&25&23\\\hline
%\textbf{2/5}&30&28&27&25&23\\\hline
%\textbf{2/6}&81&62&48&39&35\\\hline
\multicolumn{6}{|c|}{{\bf ************** Totals**************}}\\\hline
\textbf{Total Iter}&5380&1764&5287&5254&5228\\\hline
\textbf{SAM Time (s)}&1.95&3.13&6.46&12.37&25.64\\\hline
\textbf{GMRES Time (s)} &8.57&3.14&10.79&12.48&14.60\\\hline
\textbf{Total Time (s) }&10.52&6.27&17.25&24.86&40.24\\\hline
\textbf{nnz($\MN_k$)/n}&5.30&13.64&26.26&43.36&65.39\\\hline\hline
\multicolumn{6}{|c|}{{\bf ************ Totals without 2/1 ************}}\\\hline
\textbf{Total Iter}&379&327&286&253&227\\\hline
\textbf{SAM Time (s)}&1.75&2.83&5.82&11.14&23.11\\\hline
\textbf{GMRES Time (s)} &0.80&0.69&0.65&0.61&0.58\\\hline
\textbf{Total Time (s)}&2.55&3.52&6.57&11.75&23.69\\\hline
\end{tabular}
}
\caption{GMRES iterations and total runtimes for selected shifts of the Flow matrices using a priori \textit{sparsified} sparsity patterns.  A global threshold of $10^{-2}$ is used to sparsify.  An ILUTP preconditioner is computed for system 1/2, with SAM updates computed for all remaining shifts 1/3-2/6 ``Sp Patt" indicates ``Sparsified Pattern". ``SAM time" is the total amount of time spent computing all maps (but not including the time to compute the initial ILUTP).  The initial ILUTP takes 0.8 s to compute.  Totals are given for the entire sequence, as well as the sequence omitting system 2/1. (5001) indicates GMRES did not converge in the maximum allowed iterations.}
\label{table:ItersSparse}
\vspace{-7mm}
%}
\end{table}

The data show that for subsequent shifts of increasing magnitude\footnote{Note that the relative residual norm does drop at matrix 2/2, and begins to increase again for 2/3-2/6. As IRKA iterates, the corresponding shifts from one batch to the next do not drastically change (see Table \ref{table:FlowShifts}).}, the relative residual norm grows and the recycled preconditioner $\MN_k\MP_0$ becomes less effective, resulting in more
GMRES iterations as expected.  Nevertheless, recycling preconditioners using SAMs keeps
the iteration counts relatively low
for most shifts; see Section \ref{sec:ModRed} for more details and for comparison with
the results listed here.
The data also show that, as the relative residual norm decreases for denser patterns, the number of iterations decreases as well, with the exception of batch/shift 2/1, which corresponds to a very hard system. In Tables \ref{table:Iters} and \ref{table:ItersSparse}, we also show  the cumulative results for both the entire sequence, as well as the sequence omitting batch/shift 2/1.  In these tables, we specifically show iterations for batch/shift 1/2 through 2/1, and then aggregate numbers for batch/shift 2/2-2/6.  For 2/2-2/6, the iterations for each shift (and the rate of growth of those iterations) are similar to those of 1/2-2/6.

However, for the sparsity patterns derived from higher powers of $\MA_0$,
the decrease in runtime from reduced GMRES iterations does not
outweigh the increased costs of computing more expensive SAMs.
For the Flow matrices, using the pattern of $\MA_0$ comes out
as the most efficient for total runtime (both with and without batch/shift 2/1), and we use this
pattern for the experiments in Section \ref{sec:ModRed}.  Though, if we omit batch/shift 2/1, the sparsified pattern of $\MA_0$
leads to a slightly lower overall total runtime.

Next, we analyze sparsity patterns derived from the finite element mesh from which
the matrices are derived, and we focus on subsets of the sparsity pattern of the system matrix that are much sparser than the matrix itself.   This leads to maps that are cheap to compute compared with ILU-type preconditioners, which is particularly relevant for matrices that have many nonzeros per column. For this reason, we examine a sequence of matrices arising in topology optimization \cite{Shun_PhD, WangStu_2007} that
result from discretization of the 3D linear elasticity equations on
a $100\times20\times20$ trilinear (B8) element mesh,
with three unknowns per node, $u$, $v$, and $w$, giving the displacements
in the x-, y-, and z-direction.
We order nodes and variables per node lexicographically.
The size of these matrices is $n = 132$ $300$,
a typical column has up to 81 nonzeros, and the average number of
nonzeros per column varies but is approximately $70$ (after the first few iterations). In the numerical results given in Section \ref{sec:topopt}, we also consider a second, larger mesh for topology optimization.

\begin{figure}[h!]
\captionsetup{font=scriptsize}
      \begin{subfigure}{.44\textwidth}
      \captionsetup{font=small}
      \centering
          \includegraphics[scale=0.25]{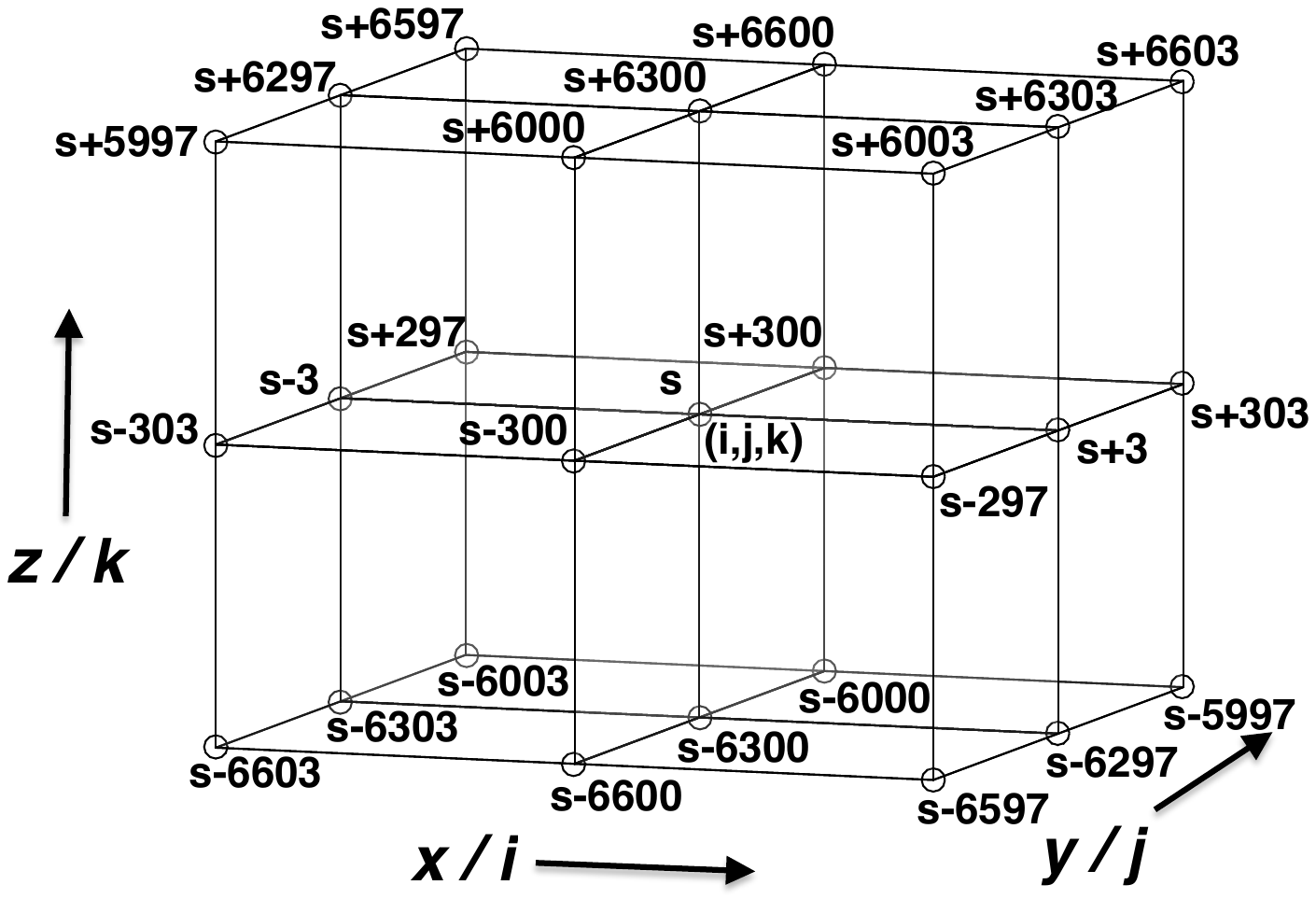}
          \subcaption{The $8$ elements containing node $(i,j,k)$ and the column indices, $s+m$, in the stiffness matrix
          that correspond to the $u$ variables of their nodes, where
          $s$ is the column index corresponding to the $u$ variable for node $(i,j,k)$.}
      \label{fig:xDisp}
      \end{subfigure}
            \hspace{8mm}
      \begin{subfigure}[t]{.44\textwidth}
      \captionsetup{font=small}
      \vspace{-3.2cm}
        \includegraphics[scale=0.35]{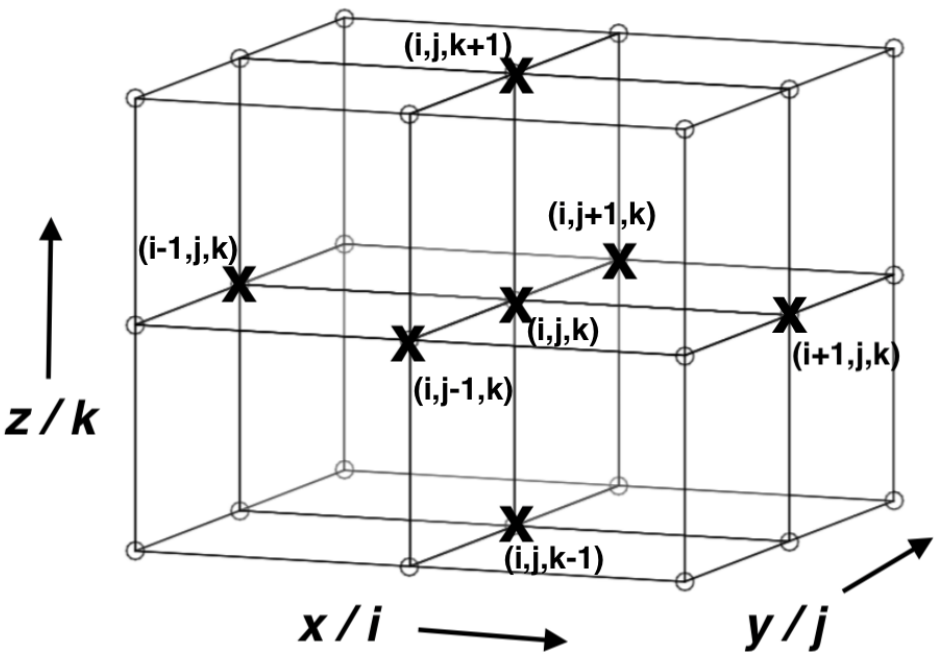}
        \vspace{1.5mm}
        \subcaption{Mesh nodes involved in the sparsity pattern Patt-0. The details are given in the text.}
        \label{fig:patt0}
    \end{subfigure}
    \vspace{-7mm}
          \caption{}
          \vspace{-5mm}
    \end{figure}

To describe the stiffness matrix derived from the mesh, we consider a typical
node $(i,j,k)$ that is not on the boundary. Let $s$, $s+1$, and $s+2$ be the
column indices corresponding to the $u$, $v$, and $w$ variables, respectively,
for node $(i,j,k)$.  Then columns $s+3$, $s+4$, and $s+5$ correspond to
the $u$, $v$, and $w$ variables, respectively, for node $(i+1,j,k)$, and columns $s-3$, $s-2$,
and $s-1$ correspond to the $u$, $v$, and $w$ variables for node $(i-1,j,k)$.
Column $s+300$ corresponds to the $u$ variable for node $(i,j+1,k)$,
and $s+6300$ corresponds to the $u$ variable for node $(i,j,k+1)$.
For the remainder of the stiffness matrix, we refer to Figure \ref{fig:xDisp},
which shows the column indices, $s+m$, corresponding to
the $u$ variables of nodes in elements that contain node $(i,j,k)$.
The column indices for the $v$ and $w$ variables at those nodes are given by
$s+m+1$ and $s+m+2$, respectively, for each $m$.

With 27 nodes, each with three displacements, there are numerous choices in sub-patterns of the matrix sparsity pattern that can be made depending on user preference.   A user may include certain displacements and omit others, dictating how many and which nonzeros will be included in the sparsity pattern.  We evaluate five of such possible patterns, Patt-0 to Patt-4. For each, we choose a selection of mesh nodes relative to $(i,j,k)$ and
displacements associated with those nodes. We stress that each of these patterns (and any of the other resulting patterns from choosing a different subset of displacements) are much sparser than the stiffness matrix itself.

To define Patt-0 for the $u$ column of the node $(i,j,k)$, corresponding to column $s$, we combine its index with the index for the $u$ variable at the nodes marked by `{\Large \textbf{$\times$}}' in Figure \ref{fig:patt0}.  The resulting sparsity pattern contains, for column $s$, the ordered pairs $(s,s)$, $(s \pm 3,s)$, $(s \pm 300,s)$, and $(s \pm 6300,s)$. For the
$v$ and $w$ columns, we use the same pattern, but with $s$
indicating the $v$, respectively, $w$ column for node $(i,j,k)$.
On boundaries, this pattern is adjusted to take the boundary and boundary conditions into account. This will also be done for Patt-1 discussed below.

Patt-1 was obtained by experimenting with minor variations of Patt-0. This sparsity
pattern contains, for column $s$, the ordered pairs $(s,s)$,
$(s \pm 1,s)$, $(s \pm 300,s)$, and $(s \pm 6000,s)$. For the
$v$ and $w$ columns, we use the same pattern, but with $s$
indicating the $v$, respectively, $w$ column for node $(i,j,k)$.  Patterns Patt-2, Patt-3, and Patt-4 include more column indices from the original sparsity pattern of the matrix based on the finite element mesh described above.\footnote{We omit specific details for Patt-2--Patt-4 here for brevity.  We refer the interested reader to \cite{GriDesGug2020} for discussions of these patterns.} All patterns, Patt-0--Patt-4, contain {\it substantially fewer} nonzeros than the stiffness matrix.  Table \ref{table:topoptPatts} shows the numbers of nonzeros in the sparsity patterns with timing and iterations results.

To evaluate the effectiveness of these patterns, we compute an ILUTP preconditioner for the matrix at optimization step 100, and
compute SAMs for the matrices at steps 105, 110, 115, 120, 125, and 130.
We solve the preconditioned systems using full GMRES.
The results are shown in Table \ref{table:topoptPatts}.
While these maps have substantially fewer nonzeros than the matrices
themselves, recycling the initial preconditioner using these SAMs
keeps the GMRES iterations low. When computing the SAM with either Patt-0 or Patt-1, each map takes less than three seconds to compute.
Including more nonzeros from the sparsity
pattern of the stiffness matrix (Patt-2--Patt-4), decreases the total
number of GMRES iterations a bit further; however, there is a substantial
increase in the time to compute the maps.  Both times and iteration counts are comparable between Patt-0 and Patt-1 as shown in Table \ref{table:topoptPatts}, though the latter is slightly faster across the longer sequence of topology optimization matrices considered in Section \ref{sec:topopt}.  So, we show results using Patt-1 for those experiments. In particular, we demonstrate the effectiveness of recycling preconditioners using SAMs in detail in Section \ref{sec:results},
providing comparisons with recomputing preconditioners and reusing preconditioners for several applications.

%\begin{center}
\begin{table}%[!ht]
\parbox{0.99\linewidth}{
\captionsetup{font=scriptsize}
%\begin{center}
\centering
\scalebox{0.7}{
\begin{tabular}{|l|c|c|c|c|c|c|}
\hline
\textbf{Optimization}&\textbf{Patt-0}&\textbf{Patt-1}& \textbf{Patt-2}&\textbf{Patt-3}&\textbf{Patt-4}\\
\textbf{Step}&\textbf{Iter}&\textbf{Iter}&\textbf{Iter}&\textbf{Iter}&\textbf{Iter}\\\hline
\textbf{100}&166&166&166&166&166\\\hline
\textbf{105}&173&173&173&170&169\\\hline
\textbf{110}&184&184&183&177&176\\\hline
\textbf{115}&195&196&194&184&183\\\hline
\textbf{120}&208&208&205&190&189\\\hline
\textbf{125}&216&220&213&194&193\\\hline
\textbf{130}&228&226&222&197&197\\\hline\hline
\textbf{Total Iter}&1370&1373&1356&1278&1273\\\hline\hline
\textbf{SAM Time (s)} &14.96&15.10 &45.47 &91.52 &  97.23 \\\hline
\textbf{GMRES Time (s)} & 97.21 &98.49 &98.19 &95.03 & 94.96 \\\hline
\textbf{Total Time (s)}&234.57&236.00 &266.06 &308.96 &314.61 \\\hline\hline
\textbf{nnz($\MN_k$)/n}&5.47&5.33&9.46&29.07&33.67\\\hline
\end{tabular}
}
\caption{GMRES iterations and total runtimes for matrices from selected steps
of the topology optimization application
using a priori mesh-based sparsity patterns.
``SAM time" is the total amount of time spent computing all maps (but not including the time to compute the initial ILUTP).
The initial ILUTP takes $122.41$s to compute.}
\label{table:topoptPatts}
}
%\end{center}
\vspace{-7mm}
\end{table}

\section{Implementation}\label{sec:impl}
We have developed efficient implementations for computing SAMs and ILUTP preconditioners as MATLAB\rtm~ m-files to make useful runtime comparisons between recycling a preconditioner and computing a new one.
We use a MATLAB\rtm~implementation of the AMG preconditioner developed as part of work done in \cite{HuLinZik2019}.  We also refer the reader to \cite{RugStu1985, RugStu1986, Stu1983, Stu2001, StuAppend2000} \cite[Appendix A: An Introduction to Algebraic Multigrid, K. St\"{u}ben]{TroOosShu2001}.%\AGM{Both Trottenberg AND Stuben's Appendix?  Or just the latter?}

First, we describe an efficient implementation for computing SAMs.
To efficiently compute the SAM updates, the solution of (\ref{eq:introSAM})
should be implemented in sparse-sparse fashion.
For most problems, the nonzero pattern of the matrices does
not change, and we preprocess the sparsity pattern for the maps
to set up data structures for the small least squares problems just once;
see Algorithm \ref{alg:SAMpre}.
Since we store $\MA_k$ as a MATLAB\rtm~  sparse matrix, access to columns is cheap.

Given a pattern, $S$, $s_\ell = \{i \,|\, (i,\ell) \in S\}$,
the set of indices of potential nonzeros in $\Vn_\ell$, column  $\ell$ of $\MN_k$.
To compute $\Vn_\ell$, we only need the $m$
columns $\Va_j$ of $\MA_k$ such that $j \in s_\ell$, and
as these columns are sparse we need only consider rows $i$ such that
$a_{i,j} \neq 0$ for some $j \in s_\ell$.
Let $r_\ell$ be the set of relevant row indices.
Then the least squares problem for $\Vn_\ell$ is defined as
\eqsn
  \Vn_\ell(s_\ell) & = & \arg \min_{\wt{\Vn} \in \Cn{m} }
    \|\MA_k(r_\ell,s_\ell) \wt{\Vn} - \MA_0(r_\ell,\ell) \|_2 ,
\eqen
where $\MA_k(r_\ell,s_\ell)$ is the submatrix of $\MA_k$ indexed
by $r_\ell \times s_\ell$,
and $\MA_0(r_\ell,\ell)$ is the corresponding
part of the $\ell$th column of $\MA_0$.  Note that if
$(a_0)_{i,\ell} \neq 0$ but $a_{i,j} = 0$
for all $j \in s_\ell$, row $i$ is irrelevant for
computing $\Vn_\ell$ since $\Ve_i \perp
\Sp{\{\Va_j \,|\, j \in s_\ell\}}$. However, if we wish
to compute, in addition to $\MN_k$, the residual
$\MRr_k = \MA_k \MN_k - \MA_0$ or its norm on Line \ref{line:Res} of
Algorithm \ref{alg:SAM}, we need to include such rows as well.
If the matrices $\MA_k$ and $\MA_0$ have the same sparsity
pattern and the pattern of $\MN_k$ includes at least the
diagonal, this is not an issue.
So, the size of each least squares problem depends only on the chosen sparsity pattern of $\MN$ and the sparsity pattern of $\MA_k$, not on the matrix size. So the least
squares problems are small, independent of $n$, and most are about the same size.
For example, for the Flow application discussed in Section \ref{sec:ModRed}, the maximum size of the least squares problems to compute the SAM updates is $20 \times 7$, independent of $n = 9~669$.

Finally, to ensure the map is computed and then stored as efficiently as
possible, in Lines \ref{line:cooStart}-\ref{line:cooEnd} of
Algorithm \ref{alg:SAM}, we compute $\MN$ in coordinate format (COO).
After the entire map has been computed, we convert this temporary data structure into a MATLAB\rtm~
sparse matrix using the command {\tt sparse} in Line \ref{line:Sparse}.

\begin{algorithm}
  \caption{Preprocessing for Computing Sparse Approximate Maps }
	\label{alg:SAMpre}
  \begin{algorithmic}[1]
  \State Given sparsity pattern $S$ and matrix $\MA$
  \State $maxSk = 0$; $maxRk = 0$; \Comment{initialize max num of columns, max num of rows}
  \For{$k = 1:n$} \Comment{for each column do}
	\State $s_k = \{i \,|\, (i,k) \in S\}$ \Comment{get indices; typically defined in advance}
    \State $r_k = \emptyset$ \Comment{Initialize set of rows for $k$th LS problem}
    \ForAll {$j \in s_k$}
      \State $t = \texttt{find}(\Va_j)$ \Comment{find indices of nonzeros in column $\Va_j$} \label{line:Find}
	  \State $r_k = r_k \cup t$
    \EndFor
    \State $nnz_k = \#( s_k )$ \Comment{\#() gives number of elements in a set}

    \hspace{-3mm} \textbf{if} {$nnz_k > maxSk$} \textbf{then} $maxSk = nnz_k$ \textbf{end if}

    \hspace{-3mm} \textbf{if} {$\#( r_k ) > maxRk$} \textbf{then}  $maxRk = \#( r_k )$ \textbf{end if}
  \EndFor
  \State Allocate $maxRk \times maxSk$ array for storing the LS matrices,
    $maxRk$ vector for storing the right hand side, and $maxSk$ vector for
    storing the solution.
  \end{algorithmic}
\end{algorithm}

\begin{algorithm}
  \caption{Computing $\MN = \arg \min_{\wt{\MN} \in \cS}\Fnm{\MA\wt{\MN} - \wh{\MA}}$}
	\label{alg:SAM}
  \begin{algorithmic}[1]
  \State $cnt = 0$ \Comment{counts number of nonzeros in preconditioner}
   \State (Preallocate space for $\MA_\mathrm{tmp}$)
  \For{$k = 1:n$} \label{line:cooStart}
	\State $\MA_\mathrm{tmp} = \MA(r_k,s_k)$ \Comment{get submatrix indexed by $r_k$ and $s_k$ for LS problem}
	\State $\Vf = \wh{\MA}(r_k,k)$ \Comment{get rhs for LS problem}
    \State Solve LS $\MA_\textrm{tmp} \Vz = \Vf$ \label{line:LS}
	\State (possibly save residual, norm of residual, etc.) \label{line:Res}
    \State $rowN[cnt+1:cnt+nnz_k] = s_k$ \Comment{assign indices in order of row ind. in $s_k$}
	\State $colN[cnt+1:cnt+nnz_k] = k$
	\State $valN[cnt+1:cnt+nnz_k] = \Vz$
	\EndFor \label{line:cooEnd}
	\State $\MN = $ {\tt sparse}$(rowN,colN,valN)$ \Comment{convert into sparse matrix} \label{line:Sparse}
  \end{algorithmic}
\end{algorithm}

Our implementation of ILUTP closely follows \cite{Saad_ILUT94,Saad09}. To make the 
implementation efficient in MATLAB\rtm, we made the following main changes. 
(1) We transpose $\MA$, in MATLAB\rtm~ sparse matrix storage, to access the rows efficiently.
(2) Where possible, we use MATLAB\rtm~ routines, such as {\tt find}, {\tt sort}, {\tt sparse}, and {\tt min}.
(3) Where possible, we have vectorized loops. 
(4) We use {\tt sparse} to build $\ML$ and $\MU$ efficiently, and we use {\tt tril} and {\tt triu} to ensure
MATLAB\rtm~ recognizes and uses the $\ML$ and $\MU$ factors as triangular matrices.
The m-file is available from \cite{ilutp-m}.

\section{Numerical Experiments} \label{sec:results}
We analyze the effectiveness of reusing, recomputing, and recycling preconditioners
for several applications. All systems are solved using preconditioned GMRES.
We compare the results of computing a new ILUTP or AMG preconditioner for each
matrix  or selected matrices, reusing the initial $\MP_0$ for all systems, updating $\MP_0$ with a
new SAM update for all systems, and updating $\MP_0$ with a SAM update only at selected
systems in the sequence, combined with computing a new $\MP_0$ at selected systems. Computing a new preconditioner for every matrix is always
the most expensive in runtime, but it provides a useful benchmark in terms of the
number of iterations. Computing a SAM update only at selected
systems (for a long sequence, combined with recomputing the preconditioner at selected systems) is usually the winner in runtime.
We report runtimes for computing
preconditioners and SAMs, the number of iterations
and runtime to solve each system, and
total runtime and number of iterations for the whole sequence.
Our focus is total runtime.

%Finally, for longer sequences of systems and other problems, many other variations of
%computing preconditioners and updates may be effective.
%\AGM{See commented code in tex file for how the following was reworded.} 
 We tested several indicators for computing a new SAM update or new preconditioner.  A simple and effective strategy is to compute a new SAM or preconditioner based on the estimated time for this computation
and the (relative) increase in the solution time for a system
or the number of iterations. For the smaller the topology optimization problem, we give results for recomputing based on a percent increase in iterations.% to the topology optimization matrices in Section \ref{sec:topopt}.}  %This strategy is a reasonable choice and we provide good results when using it.  Other strategies for recomputation may be considered, but further analysis to determine a single best indicator is left for future research.}
\comment{We have tested several indicators for computing
a new SAM update or a new preconditioner. While some results were
encouraging, we did not find a single best indicator, and we leave this for future
research.
A simple and effective strategy is to compute a new SAM or
preconditioner based on the estimated time for this computation
and the (relative) increase in the solution time for a system
or the number
of iterations.}

% REMOVED AINV - COULD RUN ON FLOW, BUT NOT CURRENTLY DONE
\comment{We also show results for the AINV preconditioner and its updates.
The AINV preconditioner was not competitive for the applications
analyzed in this paper. Hence, we show results only for the THT matrices.
Algorithms computing AINV preconditioners and AINV updates can be found in \cite{BellBert11,BenzBert03,BenzCull00,BenzHaws00,BenzKouh01,BenzTuma98,Rafi14}.}

\subsection{Topology Optimization}\label{sec:topopt}

This test problem leads to a long sequence of linear systems,
where the matrix has a relatively large number of nonzeros per row.
As a result, this problem is particularly useful to demonstrate effective
SAMs that are much sparser than the matrix itself.

Topology optimization is a structural optimization method that
optimizes the material distribution
inside a given domain~\cite{BendSig-Bk_2003, Mackerle2003, Rozvany2001,
Sigmund2000, WangStu_2007}.
The method computes a design by determining
which points of space should be material and which points
should be void (i.e. no material), combining finite element approximation,
linear solvers (for linear partial differential equations), and optimization.
In this case, we minimize the compliance (see below) subject to a volume constraint.
We specify the problem mathematically as follows.
\eqsn
  \min_{{\bf \rho},\Vu}\: &&c({\bf \rho},\Vu) = \Vu^T\MA({\bf \rho})\Vu
\label{eq:topopt}\\
  \mbox{s.t. } && \MA({\bf \rho})\Vu = \Vf \nonumber  \\
  && 0 < \rho_0 \leq \rho_e \leq 1 \qquad e=1,2,\cdots,n_e \nonumber \\
  && \int_\O \rho \, \mathrm{d} \O \leq V, \nonumber
\eqen
where $c$ is the compliance, $\MA({\bf \rho})$ is the stiffness matrix as
a function of the density distribution ${\bf \rho}$, $\Vu$ and $\Vf$ are
the displacement vector and load vector, $\rho_0$ is a chosen, small, positive
lower bound for the density to avoid singularity of the stiffness matrix, and $V$
is the total volume in use. The Solid Isotropic Material with Penalization
(SIMP) method~\cite{Bendsoe1989SIMP, Bendsoe1999SIMP} uses one design variable
to represent the density in each element.

The structure of the matrices is detailed in Section \ref{sec:theo}.
We consider matrices that result from discretizing the 3D linear elasticity equations
for variable density on a trilinear (B8) finite element mesh.  We examine  two such meshes with
 (1) $100\times20\times20$ elements,  and (2) $150 \times 30\times 30$ elements, 
with three unknowns per node, $u$, $v$, and $w$, giving the displacements
in the x-, y-, and z-directions resulting in matrices of  sizes (1) $n = 132$ $300$, and (2) $n = 432$ $450$. (see \cite{Shun_PhD, WangStu_2007} for details).  We refer to these respectively as the ``small" and ``large"  topology optimization matrices.  For both, we take a representative sequence from an optimization that typically takes hundreds, but possibly up to 1000, optimization steps before reaching the optimal design.  We demonstrate the effects of several strategies for reusing, recycling, and recomputing the preconditioner. We use the ILUTP preconditioner for the small matrices and the AMG preconditioner for the large matrices. Results are provided
in Tables \ref{table:LargeTopOpt_Reuse}-\ref{table:SAM_ILUTP_Recomp}.

 For the large matrices, we show that in the case of an expensive preconditioner resulting in fast GMRES convergence, we can preserve this good convergence with SAM updates. We use a path cover adaptive algebraic multigrid \cite{HuLinZik2019} for the preconditioner. The grid and operator complexities are 1.166 and 1.863, respectively.\footnote{The grid complexity gives the ratio of the number of all grid points to that on the finest level, and the operator complexity gives the ratio of the number of nonzero entries in all operators to that on the finest level \cite{Stu2001}\cite[p. 487, Appendix A: An Introduction to Algebraic Multigrid, K. St\"{u}ben]{TroOosShu2001}.}  
In particular,  
both are less than 2.0. %\cite{Stu2001}\cite[Appendix A: An Introduction to Algebraic Multigrid, K. St\"{u}ben]{TroOosShu2001}. 
The computation of this preconditioner takes about 16.5 minutes, whereas the SAM updates take less than five seconds, with five additional seconds for preprocessing at the first update.  We  use  full  GMRES with  maximum iterations set to 400 and the zero initial guess.  

Table \ref{table:LargeTopOpt_Reuse} shows the results when we reuse the initial AMG preconditioner. In this case, the GMRES iterations grow quickly and reach the maximum allowed twice.
At that point, we compute a new preconditioner and reuse it for subsequent systems. So, we recompute the preconditioner halfway\footnote{If we let GMRES continue beyond the maximum number of iterations for system 65, convergence is reached at 435 iterations.} through the sequence and would do this again at optimization step 70.  Total computation time is nearly two hours.
Apparently, modest changes in material distribution may already make the coarse grid correction in the AMG preconditioner substantially less effective.

Table \ref{table:LargeTopOpt_SAM} shows the results when we update the initial AMG preconditioner using SAMs at every step. We use Patt-1, described in Section \ref{sec:theo}, which has seven nonzeros in a typical column compared with up to 81 nonzeros in a typical column of the stiffness matrix. As a single GMRES iteration takes about two seconds (in both tests) due to the ILU-smoothing in the AMG preconditioner, preserving a low number of iterations is important. Recycling the initial preconditioner using SAM updates achieves this goal, with the number of iterations growing slowly from 35 to 101, avoiding recomputations of the preconditioner. The total computation time when recycling the AMG preconditioner using SAMs is reduced to less than 38 minutes.

\begin{table}
\parbox[t]{0.45\linewidth}{
\captionsetup{font=scriptsize}
%\begin{center}
\centering
\scalebox{0.7}{
\begin{tabular}{|l|c|c|c|}
\hline
\textbf{Mats}&\textbf{Prec }&\textbf{GMRES }&\textbf{Iter}\\
&\textbf{(s)}&\textbf{(s)}&\\\hline
\textbf{60}&994.52&80.34&35\\\hline
\textbf{61}&0&127.03&58\\\hline
\textbf{62}&0&222.88&99\\\hline
\textbf{63}&0&416.09&171\\\hline
\textbf{64}&0&741.18&299\\\hline
\textbf{65}&0&1026.75&400\\\hline
\textbf{66}&1002.24&83.57&37\\\hline
\textbf{67}&0&168.36&77\\\hline
\textbf{68}&0&623.67&259\\\hline
\textbf{69}&0&1026.96&400\\\hline\hline
\textbf{Total}&\multicolumn{2}{|c|}{6526.20}&1835\\\hline
\end{tabular}
}
\caption{Runtimes and iterations for the large top opt systems, when reusing the AMG preconditioner for each system after the first. When max iterations (400) is reached, we recompute the AMG preconditioner for the next system and reuse this for subsequent systems.}
\label{table:LargeTopOpt_Reuse}
}
\hspace{10mm}\parbox[t]{0.45\linewidth}{
\captionsetup{font=scriptsize}
%\begin{center}
\centering
\scalebox{0.7}{
\begin{tabular}{|l|c|c|c|}
\hline
\textbf{Mats}&\textbf{Prec }&\textbf{GMRES }&\textbf{Iter}\\
&\textbf{(s)}&\textbf{(s)}&\\\hline
\textbf{60}&994.52&80.34&35\\\hline
\textbf{61}&10.52&81.92&38\\\hline
\textbf{62}&4.97&88.08&41\\\hline
\textbf{63}&4.89&94.31&44\\\hline
\textbf{64}&4.75&101.35&47\\\hline
\textbf{65}&4.80&112.10&52\\\hline
\textbf{66}&4.89&126.10&58\\\hline
\textbf{67}&4.92&147.39&68\\\hline
\textbf{68}&4.95&177.88&81\\\hline
\textbf{69}&4.93&221.85&101\\\hline\hline
\textbf{Total}&\multicolumn{2}{|c|}{2275.47}&565\\\hline
\end{tabular}
}
\caption{Runtimes and iterations for the large top opt systems, when computing the AMG preconditioner for the first system and a SAM update for each system after the first. SAMs use the sparsity pattern Patt-1.}
\label{table:LargeTopOpt_SAM}
}
%\vspace{-3mm}
\end{table}

For the small top opt matrices, computing the ILUTP preconditioner takes about two minutes.\footnote{Using {\tt ilu} (with type `ilutp') in MATLAB\rtm~ also takes about two minutes, with approximately the same number of nonzeros in the $\ML$ and $\MU$ factors, and results in a similar number of GMRES iterations. The matrices are SPD but far from
diagonally dominant, and MATLAB\rtm's IC(0) results in a poor preconditioner, while its incomplete Cholesky with threshold fails with negative pivots.\label{foot:ILU}}
Fill in for ILUTP is set to 250 and the drop tolerance is $10^{-3}$. 
We use Patt-1 for the SAMs, as it gives slightly faster runtimes then Patt-0 for the sequence of systems considered here. 
Computing a SAM update takes a bit over two seconds. 
We use full GMRES with maximum number of iterations set to $1000$. When computing an ILUTP preconditioner for top opt matrices, diagonal scaling is essential to mitigate the huge ratio in local stiffness \cite{WangdeSt09}. This scaling varies with each matrix and combined with pivoting complicates using the previous solution as initial guess, so we use the zero initial guess.
Table \ref{table:TopOptReuse} provides the results
for reusing the initial preconditioner for all systems. Results for computing a new preconditioner for every system are not shown,
given the high cost of computing the ILUTP. Computing the SAM update for every matrix after the first reduces iterations substantially, but is still too expensive in runtime.
Therefore, we compute a SAM update for each tenth system and reuse the recycled preconditioner until the next update. Results are provided in Table \ref{table:TopOptSAM10}. 

\begin{table}
\parbox[t]{0.37\linewidth}{
\captionsetup{font=scriptsize}
%\begin{center}
\centering
\scalebox{0.7}{
\begin{tabular}{|l|c|c|c|}
\hline
\textbf{Mats}&\textbf{Prec}&\textbf{GMRES}&\textbf{Iter}\\
&\textbf{(s)}&\textbf{(s)}&\\\hline
\textbf{40-49}&120.85 &134.98 &1909\\\hline
\textbf{50-59}&0 &202.42 &2667\\\hline
\textbf{60-69}& 0&272.05 &3353\\\hline
\textbf{70-79}&0 &331.72 &3880\\\hline
\textbf{80-89}&0&396.47 &4405 \\\hline
\textbf{90-99}&0&892.64 &7405 \\\hline
\textbf{100-109}&0&1375.30 &(9779)** \\\hline
\textbf{110-119}& 0&1430.16 &(10010)*\\\hline
\textbf{120-129}&0&1432.25 &(10010)* \\\hline
\textbf{130-139}& 0&1436.56 &(10010)*\\\hline
\textbf{140-149}& 0&1429.65 &(10010)*\\\hline
\textbf{150-166}&0 &2454.73 &(10010)*\\\hline\hline
\textbf{Totals}&\multicolumn{2}{|c|}{11909.79}&90455\\\hline
\end{tabular}
}
\caption{ Runtimes and iterations for the small top opt systems with the initial ILUTP reused for all systems. Results are given in groups of ten except for (150-166). **For only one system converge reached within max iterations. *Convergence not reached for any system within max iterations. }
\label{table:TopOptReuse}
}
\hspace{8mm}
\parbox[t]{0.58\linewidth}{
\captionsetup{font=scriptsize}
%\begin{center}
\centering
\scalebox{0.7}{
\begin{tabular}{|l|c||c|c||c|c|}
\hline
\textbf{Mats}&\textbf{Prec}&\textbf{GMRES}&\textbf{Iter}&\textbf{Gain}&\textbf{Gain}\\
&\textbf{(s)}&\textbf{(s)}&&GMRES (s)&\textbf{Iter}\\\hline
\textbf{40-49}&120.85&134.98&1909&0&0\\\hline
\textbf{50-59}&2.41&189.72&2495&--12.71&--172\\\hline
\textbf{50-69}&2.19&241.57&3012&--30.48&--341\\\hline
\textbf{70-79}&2.15&283.75&3396&--47.97&--484\\\hline
\textbf{80-89}&2.13&333.22&3823&--63.25&--582\\\hline
\textbf{90-99}&2.14&361.81&4067&--530.83&--3338\\\hline
\textbf{100-109}&2.14&404.84&4414&--970.46&--5365\\\hline
\textbf{110-119}&2.14&476.00&4799&--954.16&--5211\\\hline
\textbf{210-129}&2.29&446.45&4689&--985.80&--5321\\\hline
\textbf{130-139}&2.15&465.68&4844&--970.88&--5166\\\hline
\textbf{140-149}&2.18&535.47&5324&--894.18&--4686\\\hline
\textbf{150-166}&4.30&1689.75&(13187)*&--903.72&--3830\\\hline\hline
\textbf{Totals}&\multicolumn{2}{|c|}{ 5710.32}&55959&--6225.71&--34496\\\hline
\end{tabular}
}
\caption{ Runtimes and iterations for the small top opt systems
with the initial ILUTP recycled with a SAM update for every tenth system, using Patt-1. Results are given in groups of ten except for (150-166). ``Gain" indicates the decrease in GMRES time and iterations compared with reusing the initial ILUTP. There are two SAM updates in (150-160).  *For one matrix convergence not reached within max iterations.}
\label{table:TopOptSAM10}
}
%\vspace{-5mm}
\end{table}

We can improve runtimes by recomputing the ILUTP preconditioner and the SAMs based on increases in iterations. We refer to this strategy as the ``dynamic strategy". 
Considering only the ILUTP preconditioner, we recompute the preconditioner when the number of iterations at an optimization step increases by 50\% over those from the latest ILUTP computation. Table \ref{table:ILUTP_Recomp} shows the results. In addition, we can compute a SAM update when the number of iterations at an optimization step increases by 20\% over those from the latest ILUTP computation.
We only compute one SAM update between ILUTP computations.
The results are given in Table \ref{table:SAM_ILUTP_Recomp}.   

\begin{table}
\parbox[t]{0.44\linewidth}{
\captionsetup{font=scriptsize}
%\begin{center}
\centering
\scalebox{0.7}{
\begin{tabular}{|l|c|c|c||c|}
\hline
\textbf{Mats}&\textbf{Prec}&\textbf{GMRES}&\textbf{Iter}&\textbf{ILUTP}\\
&\textbf{(s)}&\textbf{(s)}&&\textbf{at Step}\\\hline
\textbf{40-53}&120.85 &210.85 &2891 &40\\\hline
\textbf{54-74}&120.92 &320.93 &4402 &54\\\hline
\textbf{75-98}&123.24 &358.48 &4921 &75\\\hline
\textbf{99-132}&120.26 &524.04 &7181 &99\\\hline
\textbf{133-162}&121.42 &433.29 &5827 &133\\\hline
\textbf{162-166}&121.21 &46.49 &667 &162\\\hline\hline
\textbf{Totals}&\multicolumn{2}{|c|}{2621.90}&26055 &($\times 6$)\\\hline
\end{tabular}
}
\caption{Runtimes and iterations for the small top opt matrices with the dynamic strategy for ILUTP only. Each group of systems starts with a new ILUTP, as indicated in the last column.}
\label{table:ILUTP_Recomp}
}
\hspace{3mm}
\parbox[t]{0.54\linewidth}{
\captionsetup{font=scriptsize}
%\begin{center}
\centering
\scalebox{0.7}{
\begin{tabular}{|l|c|c|c||c|c|}
\hline
\textbf{Mats}&\textbf{Prec}&\textbf{GMRES}&\textbf{Iter}&\textbf{ILUTP}&\textbf{SAM}\\
&\textbf{(s)}&\textbf{(s)}&&\textbf{at Step}&\textbf{at Step}\\\hline
\textbf{40-54}&123.29&223.513&3057&40&47\\\hline
\textbf{55-77}&122.52&358.62&4735&55&66\\\hline
\textbf{78-103}&122.65&391.04&5297&78&90\\\hline
\textbf{104-144}&122.31&610.46&8305&104&122\\\hline
\textbf{145-166}&122.95&295.21&4083&145&160\\\hline\hline
\textbf{Totals}&\multicolumn{2}{|c|}{ 2492.50}&25477&($\times 5$)&($\times 5$)\\\hline
\end{tabular}
}
\vspace{1mm}
\caption{Runtimes and iterations for the small top opt matrices with the dynamic strategy for both ILUTP and SAM updates with Patt-1. 
Each group of systems starts with a new ILUTP, and includes one SAM update, as indicated in the last two columns.}
\label{table:SAM_ILUTP_Recomp}
}
\vspace{-7mm}
\end{table}

The dynamic strategy for both ILUTP computation and SAM update reduces the total number of GMRES iterations (by 578 iterations), the overall runtime (by 129 seconds), and the number of ILUTP computations (by 1) compared with the dynamic strategy for the ILUTP preconditioner only. Over the full sequence of optimization steps, these numbers increase proportionally.
%using the dynamic strategy only for ILUTP computation leads to several more, expensive, ILUTP computations. compared with both dynamically recomputing and updating with SAMs.} 

\subsection{Interpolatory Model Reduction}\label{sec:ModRed}
Our next set of linear systems arises in
the Iterative Rational Krylov Algorithm (IRKA)  \cite{GugeAnth08,AntBG20}
for computing $H_2$-optimal interpolatory reduced order models \cite{AntBG20}.
Model reduction aims to replace a
high-dimensional linear dynamical system (here single-input/single-output)
\begin{equation} \label{eqn:H}
  \ME\dot{\Vx}(t) + \MA\Vx(t) = \Vb u(t), \hspace{4mm} y(t) = \Vc^T\Vx(t),
\end{equation}
where $\ME, \MA \in \Rmn nn$, $\Vb, \Vc \in \Rn n$, input and output $u(t), y(t) \in \Rl$,
and state vector $\Vx(t)\in \Rn n $, and $n$ is very large,
by a much lower dimensional dynamical system
\begin{equation}  \label{eqn:Hr}
  \ME_r\dot{\Vx}_r(t) + \MA_r\Vx_r(t) = \Vb_r u(t), \hspace{4mm} y_r(t) = \Vc_r^T\Vx_r(t),
\end{equation}
where $\ME_r, \MA_r \in \Rmn rr$, $\Vb_r, \Vc_r \in \Rn r$, and $r \ll n$,
such that $y_r(t) \approx y(t)$ in an appropriate norm for a wide range of input selections $u(t)$.
Here, for brevity, we consider single-input/single-output dynamical systems, i.e., $u(t)$ and $y(t)$ are scalar valued functions. Discussion similarly extends to the multi-input/multi-output case.
Dynamical systems with large state-space dimension $n$ appear in many applications,
ranging from nonlinear parameter inversion to optimal control to circuit design.
The repeated simulation of these large systems may be infeasible, but
model reduction allows us to do sufficiently accurate simulations with a much smaller system.

The reduced model quantities in (\ref{eqn:Hr}) are obtained by construction of the matrices
$\MV_r, \MW_r \in \Rmn nr$ (the model reduction bases) and a Petrov-Galerkin projection
%\begin{linenomath}
\eqs  \label{eq:NewMatsVecs}
  \MA_r = \MW_r^T\MA\MV_r, \hspace{4mm} \ME_r = \MW_r^T\ME\MV_r, \hspace{4mm}
    \Vb_r = \MW_r^T\Vb, \hspace{4mm} \Vc_r = \MV_r^T\Vc.
\eqe
%\end{linenomath}
Model reduction approaches differ in their choices of $\MV_r$ and $\MW_r$
\cite{Ant2005,AntBG20,BBF2014,BGW2015,siammorbook2017,quarteroni2015reduced,hesthaven2016certified}.
In interpolatory model reduction, $\MV_r$ and $\MW_r$ are constructed
so that the reduced model transfer function $H_r(s) = \Vc_r^T(s\ME_r-\MA_r)^{-1}\Vb_r$ 
is a rational interpolant 
of the full model
transfer function $H(s) = \Vc^T(s\ME-\MA)^{-1}\Vb$.
In this case, we focus on the case that  
$H_r(s)$ is a Hermite interpolant to $H(s)$ as this is required for optimality \cite{GugeAnth08,AntBG20}. Therefore, the goal is to 
 construct $\MV_r$ and $\MW_r$ such that 
$H_r(s_j) = H(s_j)$ and $H'_r(s_j) = H'(s_j)$ for some given set
of points $s_1, \ldots, s_r$.
Constructing $\MV_r$ and $\MW_r$ as
%\begin{linenomath}
\begin{align}\label{eq:V}
  \MV_r &= [(s_1\ME-\MA)^{-1}\Vb,\dots,(s_r\ME-\MA)^{-1}\Vb], \\
\label{eq:W}
  \MW_r &= [(s_1\ME^T-\MA^T)^{-1}\Vc,\dots,(s_r\ME^T-\MA^T)^{-1}\Vc],
\end{align}
%\end{linenomath}
%and obtain the reduced model using  (\ref{eq:NewMatsVecs}).
achieves this; see \cite{AnthBeat10,AntBG20} for more details.
IRKA \cite{GugeAnth08} finds the optimal interpolation points
by alternatingly computing (\ref{eq:V})--(\ref{eq:W}) for a given set of
interpolation points $\{s_j\}$ and computing a new set $\{s_j\}$ given $\MV_r$ and $\MW_r$;
for details and some variants of IRKA, see \cite{GugeAnth08,AntBG20,BeaG12,HokMag2018,xu2010optimal,bunse-gerstner2009hom,wyatt2012issues}.
Since IRKA may take many iterations to converge to the final set of
interpolation points $\{s_j\}$, many shifted systems must be solved in
computing (\ref{eq:V}) and (\ref{eq:W}).
Each set of shifts for an IRKA iteration is called a {\em batch}.

Efficient solution of (\ref{eq:V})--(\ref{eq:W}) is an important research topic.
In \cite{BeatGuge12}, inexact solves within a Petrov-Galerkin framework are used.
In \cite{AhujdeSt12}, the recycling BiCG algorithm is proposed for model reduction, which is extended to recycling BiCGSTAB for
parametric model reduction in \cite{AhujBenn15}.
Further discussion of recycling Krylov subspace methods for model reduction applications can
be found in  \cite{FengBenn09, FengBenn13}.

We give results for one set of matrices, Flow, from \cite{reporedirect15}.
These matrices arise in a simulation of the heat exchange between a
solid body and a fluid flow, and background on this benchmark problem can be found in \cite{reporedirect15}. Rather than using computational fluid dynamics,
which is quite expensive, a flow region with a given
velocity profile is used \cite{MoosRudn04}. However, this requires
a much larger number of elements, and model reduction is used to
make the simulation efficient \cite{MoosRudn04}.
For more information see  \cite{KorRud2005,MoosRudn04,reporedirect15, RudnKorv02}.
The matrices $\MA,\ME\in \Rmn nn$ are sparse with $n=9~669$. Although $\MA$ is not symmetric,
it turns out that the shifts remain real for the three steps of IRKA used here.  We use three batches of six shifts, which are real and range from $O(1)$ to $O(10^4)$.
The shifts for the Flow matrices are provided in Table \ref{table:FlowShifts}.

We use GMRES(200) for this application  and set the maximum number of iterations to 5000.  Fill in for ILUTP is 56 and the drop tolerance is $10^{-3}$.
The pattern of $\MA_0$ is used for the SAM updates.
Results are given in Tables \ref{table:Flow_ILUTPall}--\ref{table:Flow_ILUTPSAM_select}.

An interesting case arises here. The number of iterations for the first
preconditioned system, $\MA_0 \MP_0$, is very high. For this first shift, the matrix is ill-conditioned and not diagonally dominant, and this leads to a poor ILUTP preconditioner. It makes no sense to reuse or recycle this preconditioner.
%(see Tables \ref{table:Flow_ILUTPonceP0} and
%\ref{table:Flow_ILUTPSAMP0}).
As an alternative, we compute a new 
%ILUTP 
preconditioner, $\MP_1$, for the second system.
As this preconditioned system results in fast convergence,
we recycle $\MP_1$ with SAMs for each subsequent system (see Table \ref{table:Flow_ILUTPSAMP1_select}) and at select systems (see Table \ref{table:Flow_ILUTPSAM_select}). This leads to much better
iteration counts and lower runtimes.  When computing SAM updates for select shifts, as shown in Table \ref{table:Flow_ILUTPSAM_select}, for the other shifts we reuse $\MP_1$, {\it not} the SAM updated preconditioner from a previous step. As IRKA converges, the shifts from one batch to the next do not change very much.  If the initial preconditioner works well for certain shifts in one batch, we can reuse it for the corresponding shifts in subsequent batches (see Table \ref{table:Flow_ILUTPonceP1}).  For
comparison we also provide results with $\MP_1$ reused for all subsequent systems,
leading to sightly longer runtimes than using the SAMs (see Table \ref{table:Flow_ILUTPonceP1}).
%\footnote{{\color{blue} For Table \ref{table:Flow_ILUTPonceP1}, we note the high times for the fifth and sixth shifts of each batch as compared with that of 3/1. This was consistent across many tests.  However, even when accounting for these anomalies, SAMs computed at every and select shifts are still competitive with reusing.}}  

\begin{figure}[hh]
\captionsetup{font=scriptsize}
\centering    
\includegraphics[scale=0.25]{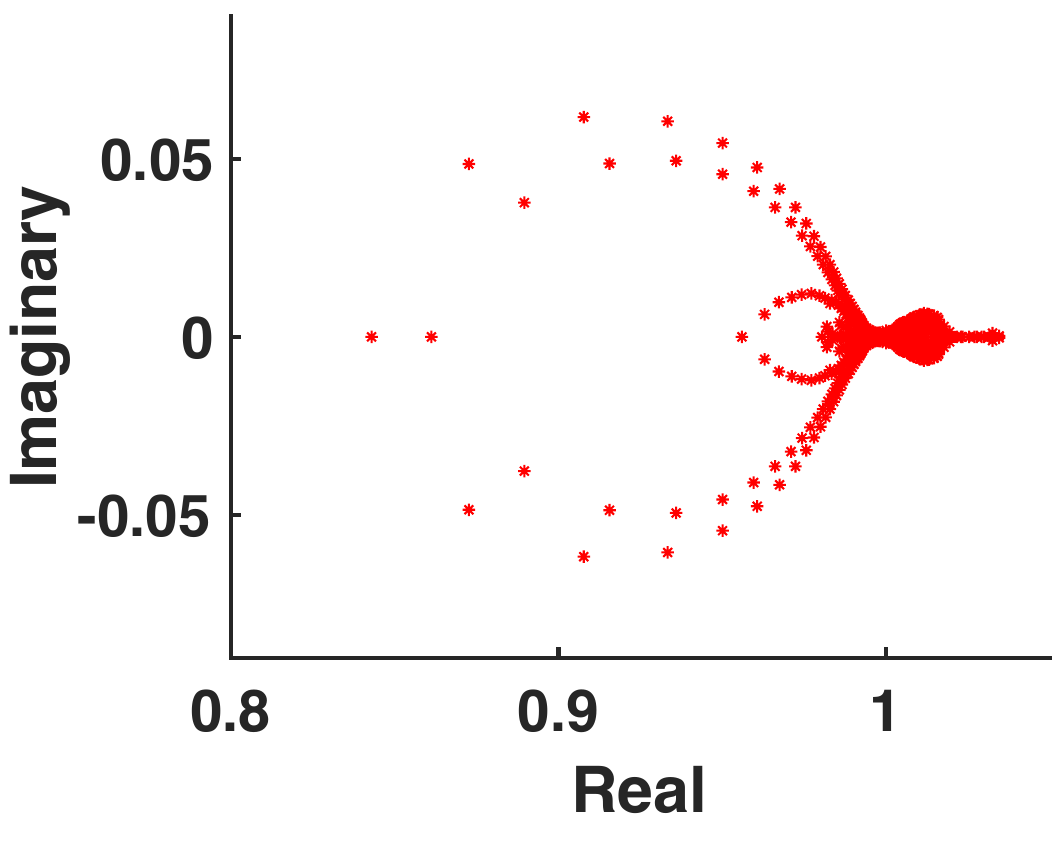}
\caption{Eigenvalues of preconditioned system 1/2, $\MA_1\MP_1$, with the preconditioner recomputed.}
\label{fig:FlowEvals12}
\end{figure}

\begin{figure}[hh]
\captionsetup{font=scriptsize}
\begin{center}  
\begin{subfigure}{.32\textwidth}
\captionsetup{font=small}	\includegraphics[width=\linewidth]{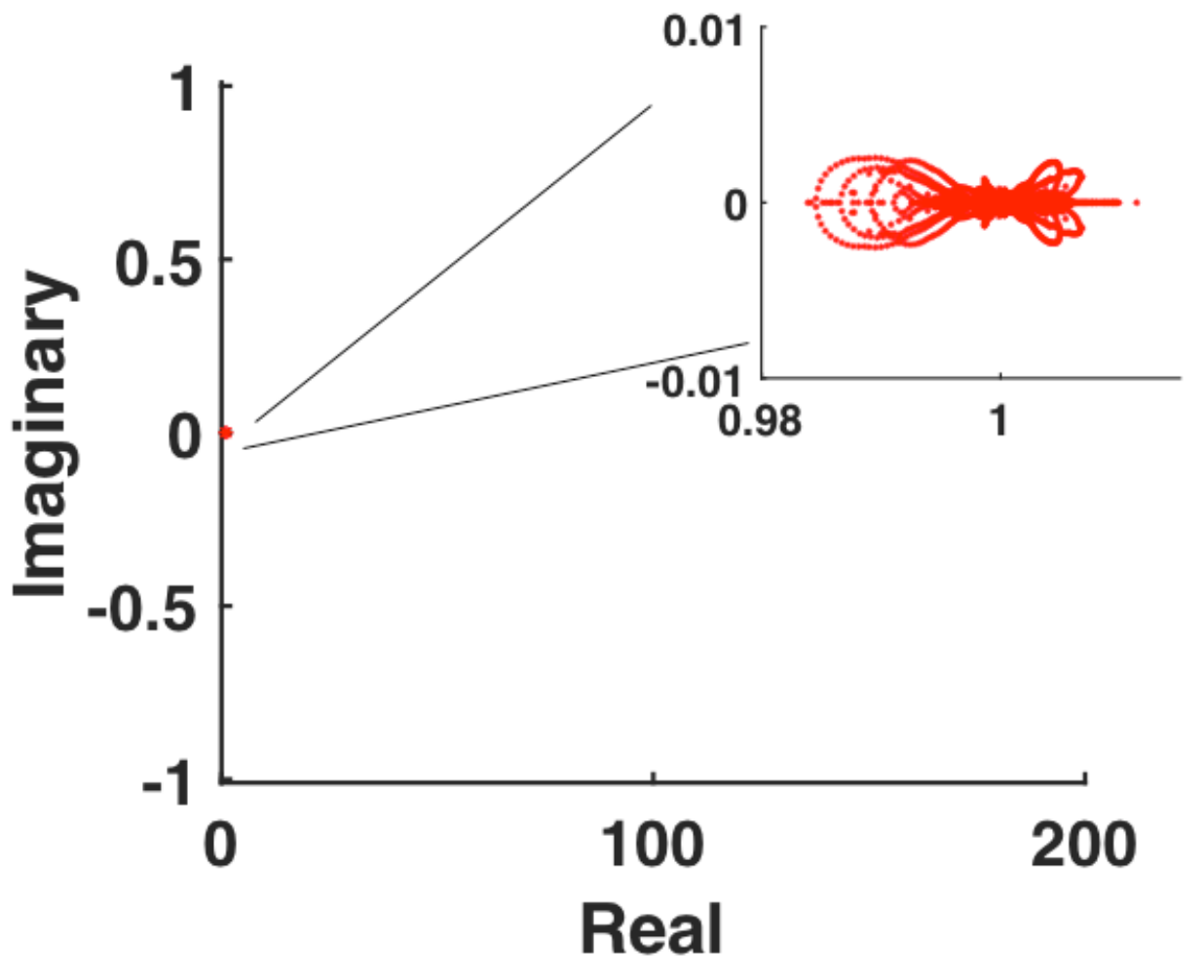}
\subcaption{}
\label{fig:Recompute5}
\end{subfigure}
\begin{subfigure}{.32\textwidth}
\captionsetup{font=small}	\includegraphics[width=\linewidth]{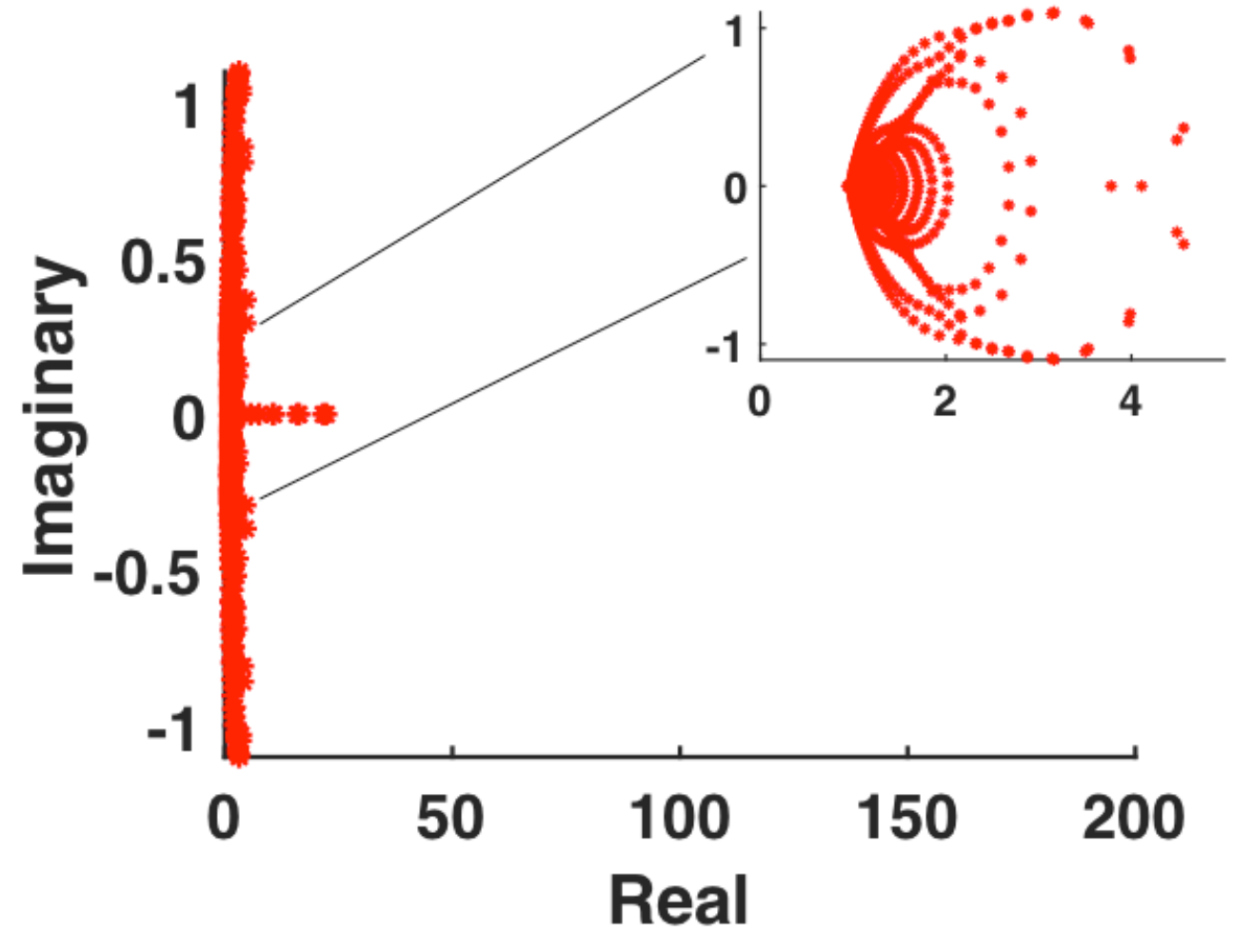}
\subcaption{}
\label{fig:Recycle5}
\end{subfigure}
\begin{subfigure}{.32\textwidth}
\captionsetup{font=small}	\includegraphics[width=\linewidth]{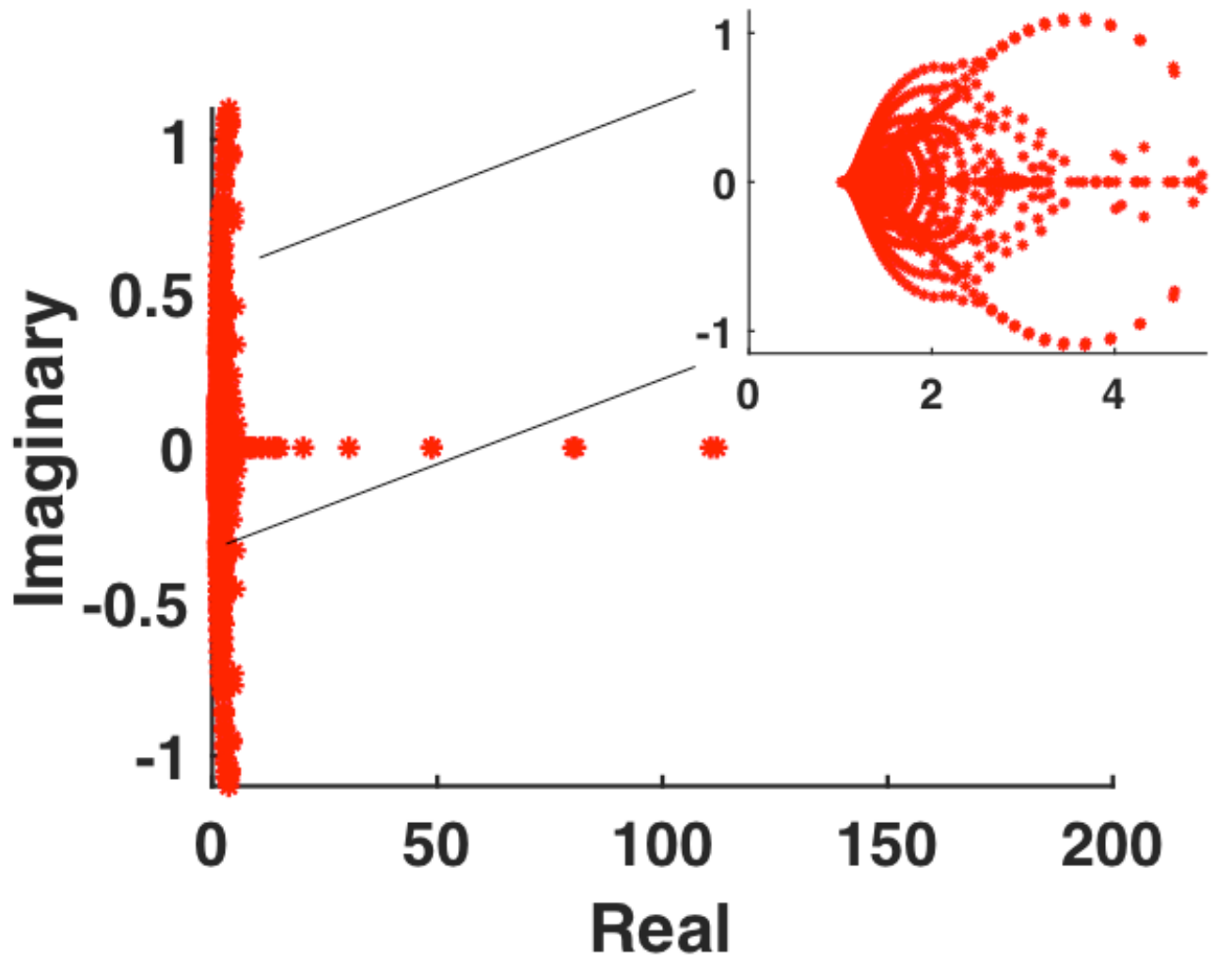}
\subcaption{}
\label{fig:Reuse5}
\end{subfigure}
\end{center}
\vspace{-5mm}
\caption{Eigenvalues of preconditioned system 1/5 with (a) $\MP_4$ recomputed, (b) $\MP_4 = \MN_4 \MP_1$ (recycled), and (c) $\MP_1$ reused. {\it Note the different scale along the real axis compared with Figure \ref{fig:FlowEvals12} and \ref{fig:FlowEvals16}.}}
\vspace{-7mm}
\label{fig:FlowEvals15}
\end{figure}

\begin{figure}[hh]
\captionsetup{font=scriptsize}
\begin{center}  
\begin{subfigure}{.32\textwidth}
\captionsetup{font=small}	\includegraphics[width=\linewidth]{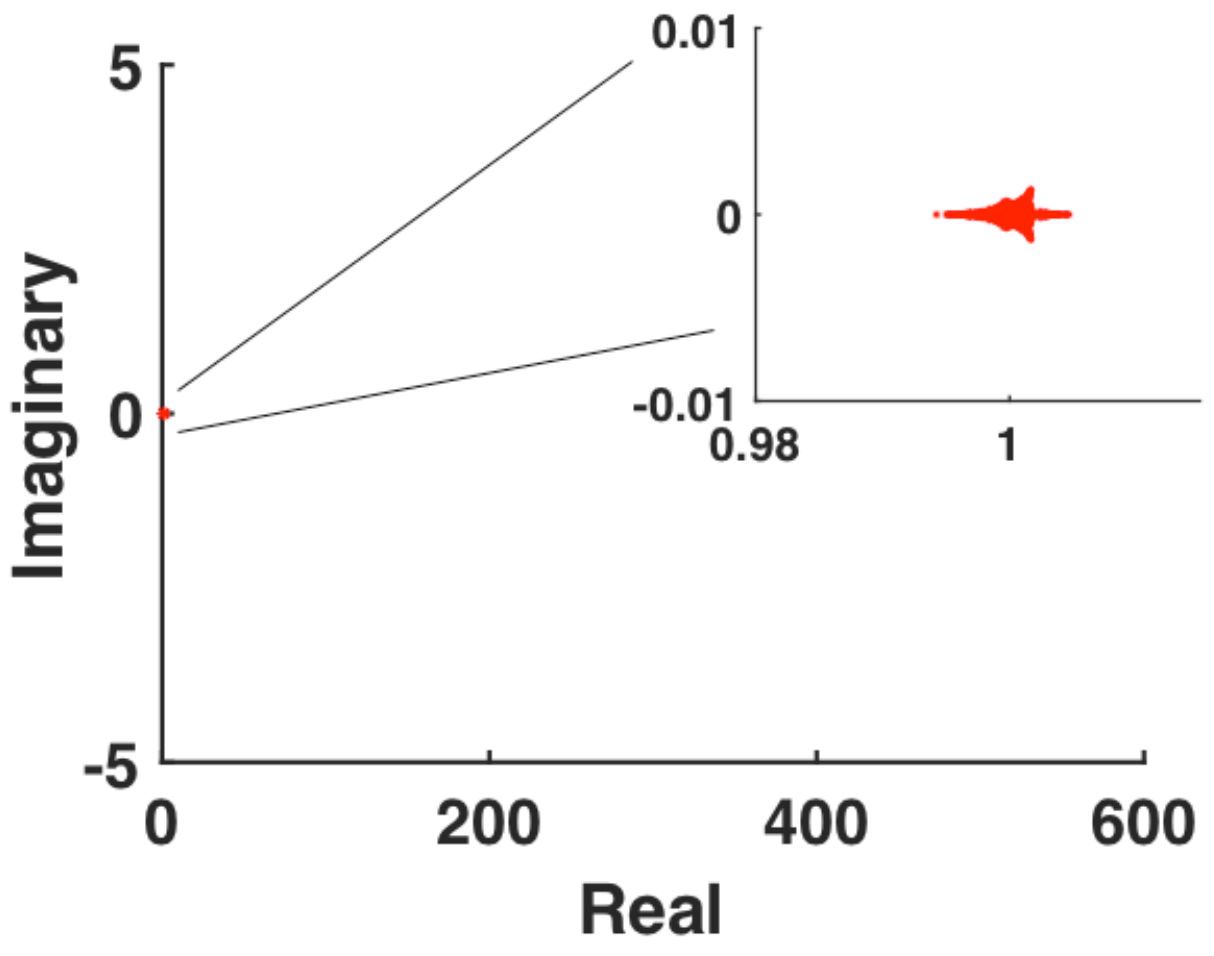}
\subcaption{}
\label{fig:Recompute6}
\end{subfigure}
\begin{subfigure}{.32\textwidth}
\captionsetup{font=small}	\includegraphics[width=\linewidth]{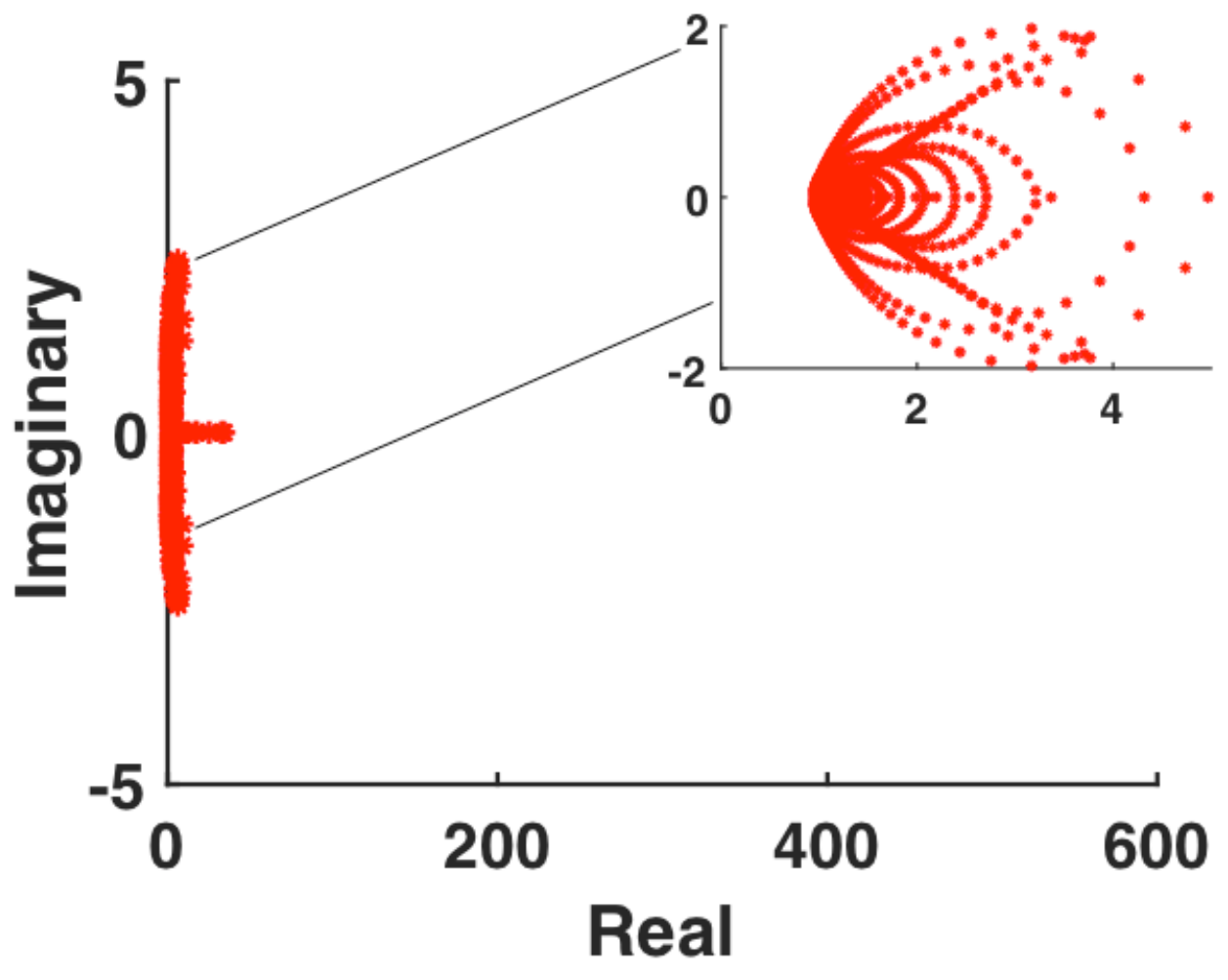}
\subcaption{}
\label{fig:Recycle6}
\end{subfigure}
\begin{subfigure}{.32\textwidth}
\captionsetup{font=small}	\includegraphics[width=\linewidth]{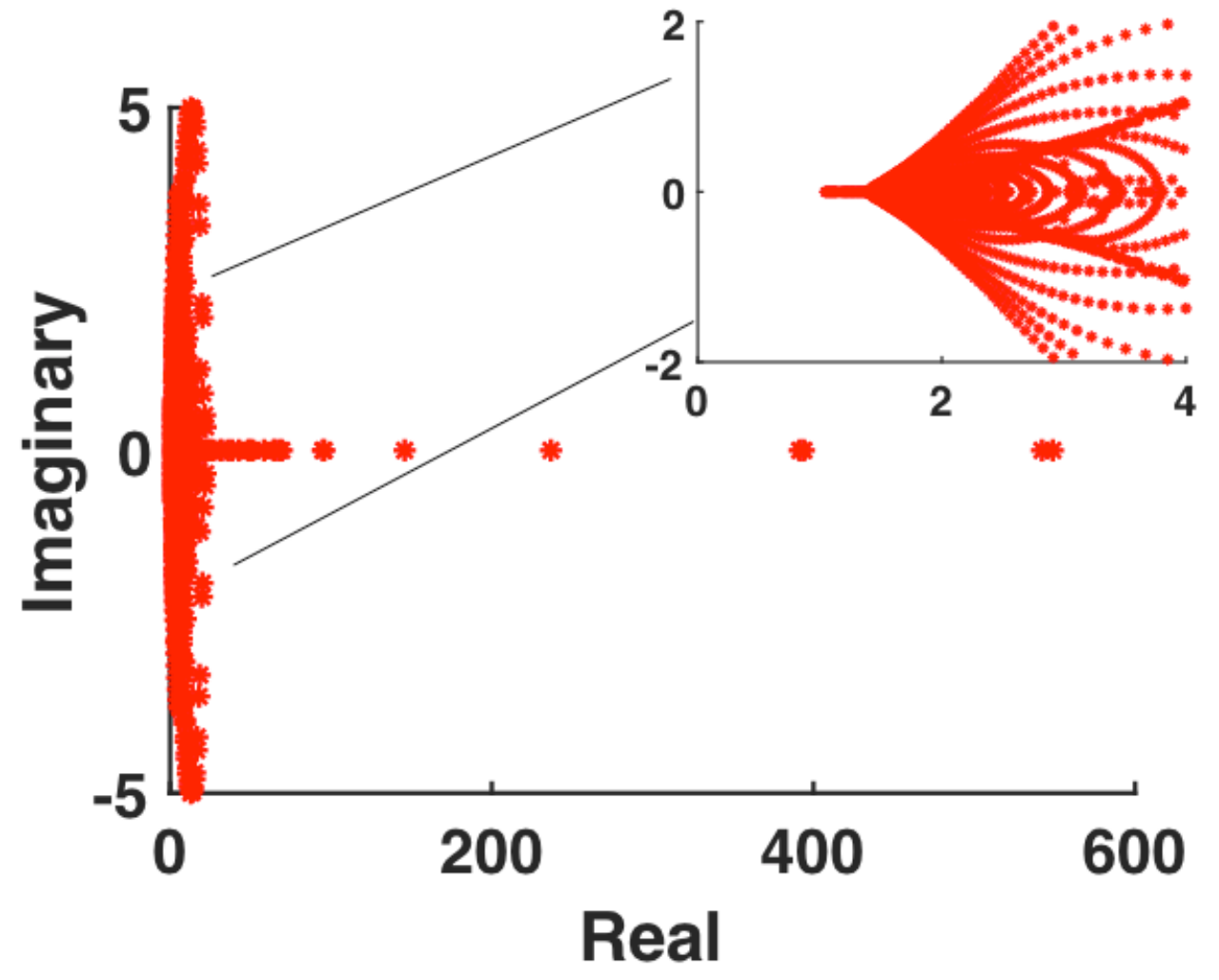}
\subcaption{}
\label{fig:Reuse6}
\end{subfigure}
\end{center}
\vspace{-5mm}
\caption{Eigenvalues of preconditioned system 1/6 with (a) $\MP_5$ recomputed, (b) $\MP_5 = \MN_5\MP_1$ (recycled), and (c) $\MP_1$ reused. {\it Note the different scale along the real axis compared with Figure \ref{fig:FlowEvals12} and \ref{fig:FlowEvals15}.}}
\vspace{-7mm}
\label{fig:FlowEvals16}
\end{figure}

As poor clustering of eigenvalues tends to lead to slow GMRES convergence, we
compare the spectra of the preconditioned systems for two shifts.   
%$\MA_k\MP_0$, where the initial preconditioner is reused, and $\MA_k\MN_k\MP_0$, where the initial preconditioner is recycled by SAMs. We also provide the spectra of $\MA_k\MP_k$ to highlight how well, compared to recomputing, both reusing and recycling preserves this clustering.  
Figure \ref{fig:FlowEvals12} shows the spectrum of the initial preconditioned system, while Figures \ref{fig:FlowEvals15} and \ref{fig:FlowEvals16} show the spectra of preconditioned systems 1/5 and 1/6, respectively, when recomputing, recycling, and reusing the preconditioner computed for system 1/2.  The figures show that the SAMs improve
the eigenvalue clustering substantially, in particular along the real axis.  At these shifts, SAMs result in much lower iterations compared with reusing the preconditioner (see Tables \ref{table:Flow_ILUTPSAMP1_select} and \ref{table:Flow_ILUTPonceP1}).

\begin{table}

\parbox{1.0\linewidth}{
\captionsetup{font=scriptsize}
%\begin{center}
\centering
\scalebox{0.7}{
\begin{tabular}{|l|c|c|c|}
\hline
\textbf{Shift }&\textbf{Batch 1}&\textbf{Batch 2}&\textbf{Batch 3}\\\hline
\textbf{1}&1.4091&1.4115&1.4116\\\hline
\textbf{2}&28.123&30.121&30.146\\\hline
\textbf{3}&150.70&163.43&163.58\\\hline
\textbf{4}&669.26&691.84&692.12\\\hline
\textbf{5}&3536.7&3565.9&3566.2\\\hline
\textbf{6}&17329&17353&17353\\\hline
\end{tabular}
}
\caption{Shifts for the Flow Matrices.  \\
For subsequent tables, B/S = Batch/Shift Number.}
\label{table:FlowShifts}
}
\parbox[t]{0.3\linewidth}{
\captionsetup{font=scriptsize}
%\begin{center}
\centering
\scalebox{0.7}{
\begin{tabular}{|l|c|c|c|}
\hline
\textbf{B/S}&\textbf{Prec }&\textbf{GMRES }&\textbf{Iter}\\
&\textbf{(s)}&\textbf{(s)}&\\\hline
\textbf{1/1}&0.80&2.64&1941\\\hline
\textbf{1/2}&0.73&0.03&13\\\hline
\textbf{1/3}&0.70&0.02&11\\\hline
\textbf{1/4}&0.67&0.02&10\\\hline
\textbf{1/5}&0.62&0.02&7\\\hline
\textbf{1/6}&0.55&0.02&7\\\hline
\textbf{2/1}&0.71&0.60&455\\\hline
\textbf{2/2}&0.70&0.03&13\\\hline
\textbf{2/3}&0.70&0.03&11\\\hline
\textbf{2/4}&0.70&0.03&11\\\hline
\textbf{2/5}&0.62&0.02&7\\\hline
\textbf{2/6}&0.55&0.02&7\\\hline
\textbf{3/1}&0.71&5.66&4231\\\hline
\textbf{3/2}&0.69&0.03&13\\\hline
\textbf{3/3}&0.73&0.02&11\\\hline
\textbf{3/4}&0.68&0.02&10\\\hline
\textbf{3/5}&0.62&0.02&7\\\hline
\textbf{3/6}&0.56&0.02&7\\\hline\hline
\textbf{Total}&\multicolumn{2}{|c|}{21.28}&6771\\\hline
\end{tabular}
}
\caption{Runtimes and iterations for the Flow matrices
with ILUTP recomputed for each shift.}
\label{table:Flow_ILUTPall}
}
\hspace{5mm}\parbox[t]{0.3\linewidth}{
\captionsetup{font=scriptsize}
%\begin{center}
\centering
\scalebox{0.7}{
\begin{tabular}{|l|c|c|c|}
\hline
\textbf{B/S}&\textbf{Prec }&\textbf{GMRES }&\textbf{Iter}\\
&\textbf{(s)}&\textbf{(s)}&\\\hline
\textbf{1/1}&0.80&2.64&1941\\\hline
\textbf{1/2}&0&0.02&18\\\hline
\textbf{1/3}&0&0.03&26\\\hline
\textbf{1/4}&0&0.59&206\\\hline
\textbf{1/5}&0&0.59&211\\\hline
\textbf{1/6}&0&0.62&232\\\hline
\textbf{2/1}&0&0.02&15\\\hline
\textbf{2/2}&0&0.02&18\\\hline
\textbf{2/3}&0&0.05&41\\\hline
\textbf{2/4}&0&0.60&207\\\hline
\textbf{2/5}&0&0.62&210\\\hline
\textbf{2/6}&0&0.78&231\\\hline
\textbf{3/1}&0&4.00&3007\\\hline
\textbf{3/2}&0&0.03&18\\\hline
\textbf{3/3}&0&0.05&40\\\hline
\textbf{3/4}&0&0.61&206\\\hline
\textbf{3/5}&0&0.63&209\\\hline
\textbf{3/6}&0&0.71&227\\\hline\hline
\textbf{Total}&\multicolumn{2}{|c|}{13.41}&7063\\\hline
\end{tabular}
}
\caption{Runtimes and iterations for the Flow matrices
with ILUTP computed for the first shift and reused for each remaining shift.}
\label{table:Flow_ILUTPonceP0}
}
\hspace{5mm}\parbox[t]{0.3\linewidth}{
\captionsetup{font=scriptsize}
%\begin{center}
\centering
\scalebox{0.7}{
\begin{tabular}{|l|c|c|c|}
\hline
\textbf{B/S}&\textbf{Prec }&\textbf{GMRES }&\textbf{Iter}\\
&\textbf{(s)}&\textbf{(s)}&\\\hline
\textbf{1/1}&0.80&2.64&1941\\\hline
\textbf{1/2}&0.23&0.03&18\\\hline
\textbf{1/3}&0.20&0.04&24\\\hline
\textbf{1/4}&0.20&0.05&30\\\hline
\textbf{1/5}&0.19&0.19&99\\\hline
\textbf{1/6}&0.19&0.14&83\\\hline
\textbf{2/1}&0.19&0.37&275\\\hline
\textbf{2/2}&0.20&0.03&18\\\hline
\textbf{2/3}&0.19&0.03&24\\\hline
\textbf{2/4}&0.19&0.25&124\\\hline
\textbf{2/5}&0.19&0.62&207\\\hline
\textbf{2/6}&0.20&0.36&155\\\hline
\textbf{3/1}&0.19&6.85&5001\\\hline
\textbf{3/2}&0.19&0.03&18\\\hline
\textbf{3/3}&0.19&0.03&24\\\hline
\textbf{3/4}&0.19&0.25&124\\\hline
\textbf{3/5}&0.19&0.19&99\\\hline
\textbf{3/6}&0.19&0.57&196\\\hline\hline
\textbf{Total}&\multicolumn{2}{|c|}{16.93}&8619\\\hline
\end{tabular}
}
\caption{Runtimes and iterations for the Flow matrices with
ILUTP computed for the first system and
SAM updates computed for all other shifts.}
\label{table:Flow_ILUTPSAMP0}
}

\parbox[t]{0.3\linewidth}{
\captionsetup{font=scriptsize}
%\begin{center}
\centering
\scalebox{0.7}{
\begin{tabular}{|l|c|c|c|}
\hline
\textbf{B/S}&\textbf{Prec }&\textbf{GMRES }&\textbf{Iter}\\
&\textbf{(s)}&\textbf{(s)}&\\\hline
\textbf{1/1}&0.80&2.64&1941\\\hline
\textbf{1/2}&0.73&0.03&13\\\hline
\textbf{1/3}&0.20&0.03&19\\\hline
\textbf{1/4}&0.20&0.05&30\\\hline
\textbf{1/5}&0.20&0.10&53\\\hline
\textbf{1/6}&0.20&0.18&79\\\hline
\textbf{2/1}&0.20&0.72&488\\\hline
\textbf{2/2}&0.20&0.02&12\\\hline
\textbf{2/3}&0.20&0.03&20\\\hline
\textbf{2/4}&0.21&0.05&30\\\hline
\textbf{2/5}&0.21&0.10&53\\\hline
\textbf{2/6}&0.19&0.18&79\\\hline
\textbf{3/1}&0.20&0.56&382\\\hline
\textbf{3/2}&0.21&0.02&12\\\hline
\textbf{3/3}&0.20&0.03&20\\\hline
\textbf{3/4}&0.20&0.05&30\\\hline
\textbf{3/5}&0.20&0.10&53\\\hline
\textbf{3/6}&0.20&0.18&79\\\hline\hline
\textbf{Total}&\multicolumn{2}{|c|}{9.83}&3393\\\hline
\end{tabular}
}
\caption{Runtimes and iterations for the Flow matrices with
ILUTP computed for the first two systems and
SAM updates computed for all other shifts.}
\label{table:Flow_ILUTPSAMP1_select}
}
\hspace{5mm}\parbox[t]{0.3\linewidth}{
\captionsetup{font=scriptsize}
%\begin{center}
\centering
\scalebox{0.7}{
\begin{tabular}{|l|c|c|c|}
\hline
\textbf{B/S}&\textbf{Prec }&\textbf{GMRES }&\textbf{Iter}\\
&\textbf{(s)}&\textbf{(s)}&\\\hline
\textbf{1/1}&0.80&2.64&1941\\\hline
\textbf{1/2}&0.73&0.03&13\\\hline
\textbf{1/3}&0&0.03&21\\\hline
\textbf{1/4}&0&0.11&60\\\hline
\textbf{1/5}&0&0.84&202\\\hline
\textbf{1/6}&0&0.89&215\\\hline
\textbf{2/1}&0&0.12&86\\\hline
\textbf{2/2}&0&0.02&12\\\hline
\textbf{2/3}&0&0.03&21\\\hline
\textbf{2/4}&0&0.06&37\\\hline
\textbf{2/5}&0&0.97&202\\\hline
\textbf{2/6}&0&0.89&215\\\hline
\textbf{3/1}&0&0.69&486\\\hline
\textbf{3/2}&0&0.02&12\\\hline
\textbf{3/3}&0&0.03&22\\\hline
\textbf{3/4}&0&0.06&36\\\hline
\textbf{3/5}&0&0.96&202\\\hline
\textbf{3/6}&0&1.01&215\\\hline\hline
\textbf{Total}&\multicolumn{2}{|c|}{10.90}&3998\\\hline
\end{tabular}
}
\caption{Runtimes and iterations for the Flow matrices
with ILUTP computed for the first two systems, and
$\MP_{1}$ reused for all remaining shifts.}
\label{table:Flow_ILUTPonceP1}
}
\hspace{5mm}\parbox[t]{0.3\linewidth}{
\captionsetup{font=scriptsize}
%\begin{center}
\centering
\scalebox{0.7}{
\begin{tabular}{|l|c|c|c|}
\hline
\textbf{B/S}&\textbf{Prec }&\textbf{GMRES }&\textbf{Iter}\\
&\textbf{(s)}&\textbf{(s)}&\\\hline
\textbf{1/1}&0.80&2.64&1941\\\hline
\textbf{1/2}&0.73&0.03&13\\\hline
\textbf{1/3}&0&0.03&21\\\hline
\textbf{1/4}&0&0.11&60\\\hline
\textbf{1/5}&0.20&0.10&53\\\hline
\textbf{1/6}&0.20&0.19&79\\\hline
\textbf{2/1}&0&0.12&86\\\hline
\textbf{2/2}&0&0.02&12\\\hline
\textbf{2/3}&0&0.03&21\\\hline
\textbf{2/4}&0&0.06&37\\\hline
\textbf{2/5}&0.21&0.10&53\\\hline
\textbf{2/6}&0.20&0.18&79\\\hline
\textbf{3/1}&0&0.68&486\\\hline
\textbf{3/2}&0&0.02&12\\\hline
\textbf{3/3}&0&0.03&22\\\hline
\textbf{3/4}&0&0.06&36\\\hline
\textbf{3/5}&0.21&0.10&53\\\hline
\textbf{3/6}&0.19&0.18&79\\\hline\hline
\textbf{Total}&\multicolumn{2}{|c|}{7.44}&3143\\\hline
\end{tabular}
}
\caption{Timings for Flow matrices with ILUTP computed for the first two systems and SAM updates computed for selected shifts.}
\label{table:Flow_ILUTPSAM_select}
}
\end{table}

\subsection{Indefinite Matrices} \label{sec:indef}
In the previous tests, computing a new ILUTP for each system gives the lowest number of GMRES iterations, but it is
too expensive in time.  Here, we consider linear
systems where the computation of the ILUTP preconditioner for the system of interest
may fail or be unstable, resulting in poor preconditioners.
This is the case, for example, for indefinite systems \cite[Chapter~10]{Saad03}.
We consider discretized 2D
Helmholtz equations $-\D u - k^2 u = f$, which arise in wave propagation problems
\cite{ErlaNabb08,ErlanVOos_04} and in flow control for unstable systems,
giving eigenvalues in both the right- and left-half planes \cite{BorgGuge14}.
In such cases, we can select 
from the set (or more generally)
a reference matrix 
for which the ILUTP algorithm computes an effective preconditioner and 
recycle this preconditioner using SAMs to
(approximately) map matrices for which ILUTP may fail to the reference matrix.

Using a modified Helmholtz equation to compute a preconditioner
has also been applied for other preconditioning approaches \cite{ErlanVOos_04}. Previous work has
successfully used operator-based preconditioners to achieve fast
convergence for Krylov methods. The shifted Laplace preconditioner \cite{ErlanVOos_04}
is used along with multilevel Krylov methods in
\cite{ErlaNabb08,ErlaVuik06,SheiLaha13}, while a  sweeping preconditioner is
constructed layer-by-layer in \cite{EngqYing11}.  Preconditioning by replacing a subset
of the Sommerfeld-type boundary conditions of the Helmholtz equation
with Dirichlet or Neumann boundary conditions is examined in \cite{ElmaOlea98,ElmaOlea99}.  We use this test problem just to demonstrate another possible use of
SAMs; we do not consider whether this approach is competitive with
the methods above.

We compute the matrix $\MK_0$ and right hand side $\Vb$ by
discretizing the 2D Laplacian on the unit square with Dirichlet boundary conditions,
$u(x,0) = 1$, $u(0,y) = 1$,  $u(x,1) = 0$, and $u(1,y) = 0$,
using a vertex-centered finite volume discretization.
$\MK_0$ is symmetric, positive definite and has size $100 \times 100$.
We compute an ILUTP preconditioner, $\MP_0$, for $\MK_0$.
Next, we solve the systems
\eqs\label{eq:helm}
  \MK_i = \MK_0-s_i\MI,
\eqe
where $\MI$ is the identity matrix, and
$s_i = i\D s$ with $\D s = 0.01$, for $i=1,2,\dots,200$.
We solve these systems with preconditioned full
GMRES, comparing a new ILUTP preconditioner
for every shift with recycling $\MP_0$ using a SAM with the pattern of $\MK_0$ for each system.
The relative convergence tolerance is $10^{-10}$ and we use a zero initial guess for each system.  For our ILUTP implementation, fill in is 20 and the drop tolerance is $10^{-3}$.  
For MATLAB\rtm's {\tt ilu} with type `ilutp', the drop tolerance is $10^{-3}$.
%The pattern of $\MK_0$ is used for the SAM updates.

\begin{figure}
\captionsetup{font=small}
\begin{center}
\begin{subfigure}{0.9\textwidth}
\captionsetup{font=small}
\includegraphics[width=\linewidth]{defToIndefBoth.pdf}
\caption{Number of GMRES iterations to converge for the discretized Helmholtz equation, comparing a new ILUTP for each $\MK_i$ (blue line) with
recycling $\MP_0$ using SAM updates for the $\MK_i$ (red line).
The results on the left are based on our ILUTP implementation. Those on the right are based on the MATLAB\rtm \ {\tt ilu} preconditioner with
type `ilutp'.}
\label{fig:indef}
\end{subfigure}
\begin{subfigure}{\textwidth}
\captionsetup{font=small}
\centering
\includegraphics[width=0.4\textwidth]{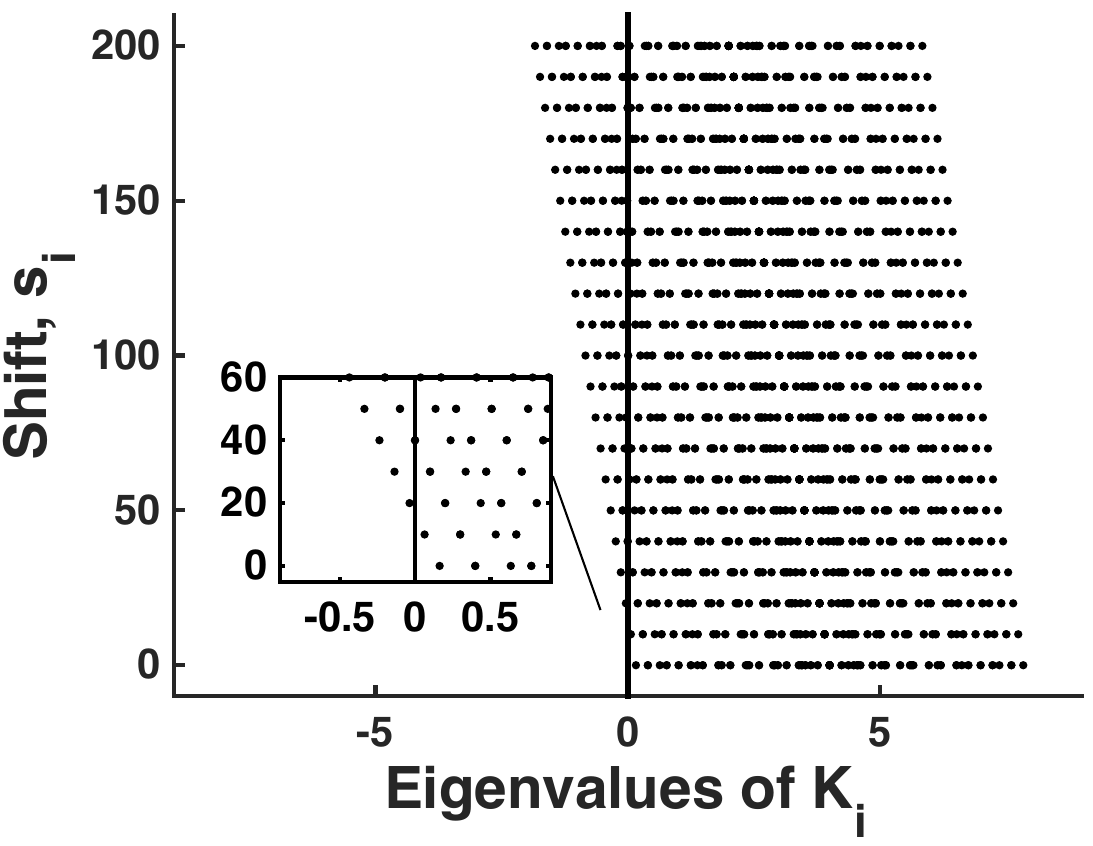}
\caption{Eigenvalues of every tenth matrix, $\MK_i$. At
the twentieth shift, the matrices become indefinite.}
\label{fig:defToIndef_eigs}
\end{subfigure}%	
\end{center}
\caption{GMRES convergence and selected eigenvalues for a discretized Helmholtz problem.}
\label{fig:helm}
\vspace{-7mm}
\end{figure}
The results are presented in Figure \ref{fig:indef}.
Figure \ref{fig:defToIndef_eigs} shows the eigenvalues for
selected $\MK_i$.
While $\MK_i$ becomes indefinite at the twentieth shift,
ILUTP produces good preconditioners until about shift $s_{125}$.
After this shift, both our
ILUTP implementation and MATLAB\rtm's $\tt ilu$ with type `ilutp'
fail to produce a good preconditioner,
and the number of GMRES iterations increases
substantially (or GMRES fails to converge).
However, using SAM updates to recycle $\MP_0$
keeps the GMRES iterations
low for almost all shifts.
For these small problems, we are not concerned with runtime and just demonstrate the superior convergence behavior obtained
with the recycled preconditioners using SAMs compared with ILUTP.

\section{Conclusions and Future Work} \label{sec:concl}
 In applications that involve many linear systems, recycling a preconditioner
can be advantageous, especially when computing a preconditioner from scratch is expensive.
We develop a flexible update to arbitrary preconditioners
that we call the Sparse Approximate Map, or SAM, update,
which can be computed for any set of closely related matrices.
The SAM is
motivated by the Sparse Approximate Inverse; however, rather
than approximately inverting a matrix, a SAM update approximately
maps a matrix to a nearby matrix for which a good preconditioner is
available.  Using SAMs, the cost of computing a very good
preconditioner can be amortized over many
systems in a sequence, since computing SAMs is cheap.
Further, a SAM is independent of preconditioner type and quality.
The sparsity patterns for SAMs can be based on powers of $\MA_0$,
mesh-based patterns, or any other salient feature of a specific problem.

In future work, we plan to consider incremental SAM updates
as in (\ref{eq:incrementSAM2}), applying maps from the left as in (\ref{eq:LeftMap}), and maps
that allow CG and MINRES to be used.
Another important future topic are approaches that
update only a few columns or rows of the map to match localized
changes in the matrix $\MA_k$. More generally, we plan to consider maps that are tuned to various structural changes in a sequence of matrices that are known a priori. We also plan to develop more indicators
for computing a new map. Finally, we plan to consider other types of maps, including different, potentially adaptive, choices in sparsity patterns.

\section{Acknowledgements} 
We thank Xiaozhe Hu for sharing his MATLAB\rtm~ implementation of the AMG preconditioner with us and for his advice how to incorporate it into our solvers.  We also thank Tania Bakhos, Arvind Saibaba, and Peter Kitanidis for providing us with matrices from a THT application \cite{Bakhos2015940}. The application is not represented in the final draft of this paper, but it was very helpful in the development of our ideas.  
We further thank the anonymous reviewers for their suggestions, which led to several valuable improvements in this paper.  
%Finally, reviews of earlier versions of this paper led to a valuable analysis of sparsity patterns, especially
%patterns that are much sparser than the matrix, as well as 
%more efficient implementations of our algorithms.}

\bibliography{paper_v02_bib}
\bibliographystyle{siam}
\end{document}